\documentclass[a4paper]{amsart}\usepackage{amsfonts,amssymb,bbm,graphicx} 
\usepackage{enumerate}

\def\C{\Bbb{C}}\def\G{\Bbb{G}}\def\k{\mathbbm{k}}
\def\bk{{\bar{\k}}}
\def\R{\Bbb{R}}\def\Z{\Bbb{Z}}

\def\li{~\\ $\bullet$ }\def\di{\partial}

\def\liml{\lim\limits}
\def\suml{\sum\limits}
\def\capl{\mathop\cap\limits}
\def\cupl{\mathop\cup\limits}
\def\wedgel{\mathop\wedge\limits}

\newcommand{\quotient}[2]{{\left.\raisebox{1.6ex}{$#1$}\!\!\!\!\!{\scalebox{2}{\ensuremath\diagup}}
\!\!\!\!\!\raisebox{-1ex}{$#2$}\right.}}
\newcommand{\quotients}[2]{{\footnotesize\left.\raisebox{0.4ex}{$#1$}\! / \!\raisebox{-0.4ex}{$#2$}\right.}}

\def\ta{\tilde{a}}\def\tA{\tilde{A}}\def\tU{{\tilde{U}}}\def\tV{\tilde{V}}
\def\tx{\tilde{x}}

\def\hA{\widehat{A}}\def\hB{\hat{B}}\def\hf{\hat{f}}\def\hG{{\widehat{G}}}\def\hI{{\hat{I}}}
\def\hR{\widehat{R}}\def\hS{\hat{S}}
\def\hphi{\widehat{\phi}}

\def\al{\alpha}\def\de{\delta}\def\De{\Delta}

\def\la{\lambda}\def\Si{\Sigma}

\def\cD{\mathcal D}
\def\cG{\mathcal G}\def\cK{{\mathcal K}}\def\cO{\mathcal O}\def\cR{\mathfrak{R}}
\def\cm{{\frak m\hspace{0.05cm}}}\def\cp{{\frak p\hspace{0.05cm}}}\def\cq{{\frak q\hspace{0.05cm}}}
\def\cU{\mathcal U}\def\cV{\mathcal V}

\def\uf{\underline{f}}\def\ug{\underline{g}}

\def\ux{\underline{x}}

\def\one{{1\hspace{-0.1cm}\rm I}}\def\zero{\mathbb{O}}

\newcommand{\ber}{\begin{array}{l}}\newcommand{\eer}{\end{array}}
\newcommand{\bpm}{\begin{pmatrix}}\newcommand{\epm}{\end{pmatrix}}
\newcommand{\bM}{\begin{matrix}}\newcommand{\eM}{\end{matrix}}
\newcommand{\bee}{\begin{enumerate}}\newcommand{\eee}{\end{enumerate}}

\def\Span{{\rm Span}}\def\Supp{{\rm Supp}}\def\iff{{\rm iff}}

\def\wrt{with respect to }
\def\sset{\subset}\def\sseteq{\subseteq}\def\ssetneq{\subsetneq}\def\smin{\setminus}

\def\Mat{Mat_{m\times n}(R)}\def\Matm{Mat_{m\times n}(\cm)}

\def\bull{\vrule height .9ex width .9ex depth -.1ex }

\newtheorem{Lemma}{Lemma}[section]\newcommand{\bel}{\begin{Lemma}}\newcommand{\eel}{\end{Lemma}}
\newtheorem{Theorem}[Lemma]{Theorem}\newcommand{\bthe}{\begin{Theorem}}\newcommand{\ethe}{\end{Theorem}}
\newtheorem{Proposition}[Lemma]{Proposition}\newcommand{\bprop}{\begin{Proposition}}\newcommand{\eprop}{\end{Proposition}}
\newtheorem{Corollary}[Lemma]{Corollary}\newcommand{\bcor}{\begin{Corollary}}\newcommand{\ecor}{\end{Corollary}}
\newtheorem{Definition}[Lemma]{Definition}\newcommand{\bed}{\begin{Definition}}\newcommand{\eed}{\end{Definition}}
\newtheorem{Definition-Proposition}[Lemma]{Definition-Proposition}

\def\bpr{~\\{\em Proof.\ }}\newcommand{\epr}{$\bull$\\}
\newtheorem{Remark}[Lemma]{Remark}\newcommand{\beR}{\begin{Remark}\rm}\newcommand{\eeR}{\end{Remark}}
\newtheorem{Example}[Lemma]{Example}\newcommand{\bex}{\begin{Example}\rm}\newcommand{\eex}{\end{Example}}
\newtheorem{Problem}[Lemma]{Problem}\newcommand{\bprob}{\begin{Problem}\rm}\newcommand{\eprob}{\end{Problem}}
\newtheorem{Properties}[Lemma]{Properties}

\newcommand{\bet}{\begin{tabular}{cccccccc}}\newcommand{\eet}{\end{tabular}}
\newcommand{\beq}{\begin{equation}}\newcommand{\eeq}{\end{equation}}

\newcommand{\bin}[2]{\binom{#1}{#2}}
\newcommand\isom{\xrightarrow{\,\smash{\raisebox{-0.65ex}{\ensuremath{\scriptstyle\sim}}}\,}}

\title[]{F\MakeLowercase{inite determinacy of matrices over local rings}.II.
\\T\MakeLowercase{angent modules to the miniversal deformations for group-actions involving the ring automorphisms}}
\author{G\MakeLowercase{enrich} B\MakeLowercase{elitskii and} D\MakeLowercase{mitry} K\MakeLowercase{erner}}
\address{Department of Mathematics, Ben Gurion University of the Negev, P.O.B. 653, Be'er Sheva 84105, Israel.}
\email{genrich@math.bgu.ac.il}
\email{dmitry.kerner@gmail.com}
\date{\today}
\thanks{D.K. was supported by the grant FP7-People-MCA-CIG, 334347. }

\begin{document}\maketitle
\begin{abstract}
We consider matrices with entries in a local ring, $\Mat$. Fix a group action, $G\circlearrowright\Mat$, and a subset of allowed deformations, $\Si\sseteq\Mat$. The standard question (along the lines of Singularity Theory) is the finite-$(\Si,G)$-determinacy of matrices.

 In 
 \cite{Belitski-Kerner}
 this determinacy question was reduced to the study of the tangent spaces $T_{(\Si,A)}$, $T_{(GA,A)}$, and their quotient, the tangent module to the miniversal deformation, $T^1_{(\Si,G,A)}=\quotients{T_{(\Si,A)}}{T_{(GA,A)}}$.
 In particular, the order of determinacy is controlled by the annihilator of this tangent module, $ann(T^1_{(\Si,G,A)})$.
 In \cite{Belitski-Kerner1} we have studied this tangent module for the group action $GL(m,R)\times GL(n,R)\circlearrowright\Mat$ and  various natural subgroups of it. These are $R$-linear group actions.

In the current work we study 
the module $T^1_{(\Si,G,A)}$ for group actions that involve the automorphisms of the ring, or,  geometrically, group-actions that involve the local coordinate changes. (These actions are not $R$-linear.)
We obtain various bounds on the support of this module. This gives ready-to-use criteria of determinacy for matrices, (embedded) modules and (skew-)symmetric forms.
\end{abstract}
\setcounter{secnumdepth}{6} \setcounter{tocdepth}{1}
\tableofcontents

\section{Introduction}
\subsection{Setup}
Let $R$ be a (commutative, associative) local ring over some base field $\k$ of zero characteristic. Denote by  $\cm\sset R$ the maximal ideal.
 (As the simplest case, one can consider
{\em regular} rings, $\cO_{(\k^p,0)}$, e.g. the rational functions regular at the origin, $\k[x_1,\dots,x_p]_{(\cm)}$ or
the formal power series, $\k[\![x_1\,\dots,x_p]\!]$. For $\k\subseteq\C$, or any other normed field,
one can consider the convergent power series, $\k\{\dots\}$ or the smooth functions $C^\infty(\R^p,0)$.)
Geometrically, $R$ is the ring of regular functions on the (algebraic/formal/analytic etc.) germ $Spec(R)$.

Let $\Mat$ be the space of $m\times n$ matrices with entries in $R$. We always assume $m\le n$, otherwise one
can transpose the matrix. Various groups act on $\Mat$.
\bex\label{Ex.Intro.Typical.Groups}
$\bullet$ The left multiplications $G_l:=GL(m,R)$, the right and the two-sided multiplications $G_r:=GL(n,R)$,  $G_{lr}:=G_l\times G_r$. Matrices considered up to $G_r$-transformations correspond to the embedded modules,  $Im(A)\sset R^{\oplus m}$. Matrices considered up to $G_{lr}$-transformations corresponds to the non-embedded modules,  $coker(A)=\quotients{R^{\oplus m}}{Im(A)}$.
\li The congruence, $G_{congr}=GL(m,R)\circlearrowright Mat_{m\times m}(R)$, acts by $A\to UAU^t$. Matrices considered up to the congruence correspond to the bilinear/symmetric/skew-symmetric forms.
\li The group of ring automorphisms, $Aut_\k(R)$,  coincides in the ``geometric case" (e.g. $R=\quotients{\k[[\ux]]}{I}$, $\quotients{\k\{\ux\}}{I}$) with the group of local coordinate changes. In singularity theory this is traditionally called ``the right equivalence", $\cR$. The automorphisms act on matrices entry-wise.
\li  Accordingly one considers the semi-direct products, $\cG_{l}:=G_l\rtimes Aut_\k(R)$, $\cG_{lr}:=G_{lr}\rtimes Aut_\k(R)$\dots Sometimes one considers only those coordinate changes that preserve some ideal $I$, i.e. the locus $V(I)\sset Spec(R)$.
\li  For a proper  ideal $J\subsetneq R$  and a group action $G\circlearrowright\Mat$ one often considers the subgroup of transformations that are trivial modulo $J$:
\beq
G^{(J)}:=\{g\in G|\ g\cdot Mat_{m\times n}(J)=Mat_{m\times n}(J)\ and\ [g]=[Id]\circlearrowright\quotient{Mat_{m\times n}(R)}{Mat_{m\times n}(J)}\}.
\eeq
For example, for $J=\cm^{k}$ the group $G^{(\cm^k)}$ consists of elements that are identities up to the order
 $(k-1)$. Similarly, $Aut^{(J)}(R)=\{\phi\in Aut_\k(R) |\ [\phi]=[Id]\circlearrowright \quotients{R}{J}\}$
\eex

\

We study deformations of matrices. In applications one often deforms a matrix not inside the whole $\Mat$ but only inside a ``deformation subspace" (a subset of prescribed deformations), $A\rightsquigarrow A+B$, $A+B\in\Si\sseteq\Mat$.

\

In this paper the subset $\Si-\{A\}\sseteq\Mat$ is a submodule. Besides the trivial choice $\Si=\Mat$, we mostly work with the spaces of symmetric matrices, $\Si=Mat^{sym}_{m\times m}(R)$, and skew-symmetric matrices, $\Si=Mat^{skew-sym}_{m\times m}(R)$.

\

 Recall that any matrix is $G_{lr}$-equivalent to a block-diagonal, $A\stackrel{G_{lr}}{\sim}\one\oplus \tA$,
where  all the entries of $\tA$ lie in the maximal ideal $\cm$, i.e. vanish at the origin of $Spec(R)$. Similar statements hold for (skew-)symmetric matrices \wrt $G_{congr}$, see e.g. \S3.1 in \cite{Belitski-Kerner1}.
 This splitting is natural in various senses and is standard in commutative algebra,  singularity theory and other fields. Often it is $\tA$ that carries the essential information.
 Therefore we often assume $A|_0=\zero$, i.e. $A\in\Matm$.

\subsection{The tangent spaces}
Fix an action $G\circlearrowright\Mat$, a deformation space $\Si\sseteq\Mat$, and a matrix $A\in\Si$. We assume that the orbit $GA$ and the deformation space possess well defined tangent spaces at $A$, which are $R$-modules.
The precise conditions/statements are in \S \ref{Sec.Background.Tangent.Spaces}.

The standard approach of deformation theory is to establish the existence of the miniversal (semi-universal) deformation and, when the later exists, to understand/to compute its tangent cone.
Accordingly one passes from the study of the germs $(GA,A)\sseteq(\Si,A)$, to the study of the tangent spaces, $T_{(GA,A)}\sseteq T_{(\Si,A)}$.

Much of the information about the deformation problem is encoded in the quotient module
\beq
T^1_{(\Si,G,A)}:=\quotient{T_{(\Si,A)}}{T_{(GA,A)}}.
\eeq
 This $R$-module is the tangent space to the miniversal deformation, when the later exists and is smooth.
In Singularity Theory such a module is known as the Tjurina algebra for the contact equivalence, and the Milnor algebra for the right equivalence.

\

This module $T^1_{(\Si,G,A)}$ is the main object of our study. We study it for those actions $G\circlearrowright\Si$ of example \ref{Ex.Intro.Typical.Groups} that involve $Aut_\k(R)$.

\subsection{(In)finite determinacy}
 The $(\Si,G)$-order of determinacy of $A$ is  the minimal number $k\le\infty$ satisfying: if $A,A_1\in\Si$ and $jet_k(A)=jet_k(A_1)$ then $A_1\in GA$. Here $jet_k$ is the projection $\Mat\stackrel{jet_k}{\to}Mat_{m\times n}(\quotients{R}{\cm^{k+1}})$. More precisely:
\bed\label{Def.Intro.Determinacy.Classical}
$ord^\Si_G(A):=min\Big\{k|\ \Si\cap\big(\{A\}+Mat_{m\times n}(\cm^{k+1})\big)\sseteq GA\Big\}\le\infty.$
\eed
\vspace{-0.3cm}(We assume that the minimum is taken here over a non-empty set.)

\

If $ord^\Si_G(A)<\infty$ then $A$ is called {\em finitely-$(\Si,G)$-determined}.
 In words: $A$ is determined (up to $G$-equivalence) by a finite number of terms in its Taylor expansion at the origin. An immediate consequence is the algebraization: $A$ is $G$-equivalent to a matrix of polynomials.

 If $ord^\Si_G(A)=\infty$ the matrix is called {\em infinitely determined}.
 This condition is empty when $\cm^\infty=\{0\}$.
 But for some rings $\cm^\infty:=\capl_{k>0}\cm^k\neq\{0\}\sset R$, then even the infinite determinacy is non-trivial  property. (The main example is the ring of germs of smooth functions, $R=C^{\infty}(\R^p,0)$, and its sub-quotients.) In this case $A$ is determined (up to $G$-equivalence) by its image under the $\cm$-adic completion, $\hA\in Mat_{m\times n}(\hR)$, i.e. its full Taylor expansion at the origin.

\

The (in)finite determinacy is the fundamental notion of Singularity Theory. More generally, as the determinacy expresses the ``minimal stability", it is important in any area dealing with matrices over rings (or matrix families or matrices depending on parameters). A trivial consequence of the finite determinacy is the algebraizability, $A$ is $G$-equivalent to a matrix of polynomials. Even more, the order of determinacy gives an upper bound on the degrees of polynomials.
 See \S\ref{Sec.Results.Relation.to.Singular.Theory} for a brief review and the relation of our work to the known results.

\subsection{Contents of the paper}

In \cite{Belitski-Kerner} we have reduced the study of determinacy to the understanding of the support/annihilator of the module $T^1_{(\Si,G,A)}$. The matrix $A$ is (in)finitely determined iff $\cm^N T_{(\Si,A)}\sseteq T_{(GA,A)}$ for some $N\le\infty$, alternatively $\cm^N T^1_{(\Si,G,A)}=\{0\}$ or $\cm^N\sseteq ann(T^1_{(\Si,G,A)})$. The order of determinacy is fixed by the annihilator $ann(T^1_{(\Si,G,A)})$, see  \S\ref{Sec.Background.BK} for more detail.

\

This is not yet the full solution of the determinacy problem, as the module $T^1_{(\Si,G,A)}$ is complicated even in the simplest cases. For example, for $A\in\Si=\Matm$,
 $T^1_{(\Si,G_{lr},A)}=Ext^1_R\Big(coker(A),coker(A)\Big)$, while the support of $T^1_{(\Si,Aut_\k(R),A)}$ is the critical locus of the map $Spec(R)\stackrel{A}{\to}Mat_{m\times n}(\k)$. The remaining question is to compute (or at least to bound) the annihilator $ann(T^1_{(\Si,G,A)})$.

\

In \cite{Belitski-Kerner1} we have studied the module $T^1_{(\Si,G,A)}$ for  group actions that do not involve $Aut_\k(R)$, in particular for
$G_r$, $G_l$, $G_{lr}$, $G_{congr}$, $G_{conj}$. We have obtained rather tight bounds on the annihilator of this module. As the immediate applications we have obtained tight bounds on the order determinacy of modules over local rings, (skew-)symmetric forms, flags of modules, flagged morphisms, chains of modules.

\

In this paper we study the module $T^1_{(\Si,G,A)}$ for the group actions $Aut_\k(R)$, $\cG_{lr}$, $\cG_{congr}$.

\

Both $T_{(\Si,A)}$ and $T_{(GA,A)}$ are (in general non-free) $R$-modules of high rank. Their quotient is usually complicated, even in the case of a regular ring, say $R=\k[[x_1,\dots,x_p]]$.
In most cases the best we can hope for is: to approximate $ann(T^1_{(\Si,G,A)})$ by some (tight) lower/upper bounds.
As the first approximation, one finds the ``set-theoretic" support of this module, i.e. the
radical of the ideal, $\sqrt{ann(T^1_{(\Si,G,A)})}$.
 This is done by checking the localizations of $T^1_{(\Si,G,A)}$ at all the possible prime ideals.
 Geometrically this amounts to checking the fibres of the sheaf, $T^1_{(\Si,G,A)}|_{pt}$, for all the points $pt\in Spec(R)$.
  The points where $T_{(GA,A)}|_{pt}\ssetneq T_{(\Si,A)}|_{pt}$  define the ``degeneracy" locus in $Spec(R)$ whose ideal is precisely the radical  $\sqrt{ann(T^1_{(\Si,G,A)})}$.

\

Section-wise the organization is as follows. The main results are stated in \S\ref{Sec.Results}.
 (We state both the bounds on the support of $T^1_{(\Si,G,A)}$ and the corresponding applications to the determinacy.)
In \S\ref{Sec.Preparations} we collect all the needed preliminaries: the algebra-geometric dictionary of localizations, the determinantal ideals and the (generalized) annihilators-of-cokernels, the singular locus of an ideal, the tangent spaces to group-orbits, the relevant approximation properties of rings.
In \S\ref{Sec.Proofs} we give all the proofs and further examples/corollaries.

\section{The main results}\label{Sec.Results}
\subsection{Notations and conventions}\label{Sec.Results.Notations}

\subsubsection{}\label{Sec.Results.Relevant.Approx}
When computing the order of determinacy (using the results of \cite{Belitski-Kerner}) we need some restrictions on the ring $R$, the so-called `{\em relevant approximation property}'. For simplicity in this paper we restrict to the following particular cases (here $\ux$ denotes a finite tuple of variables), see \S\ref{Sec.Background.Approximation.Properties} for more details.
\bee[i.]
\item Either the ring is complete, $R=\quotients{\k[[\ux]]}{I}$;
\item or $R$ is analytic, $R=\quotients{\k\{\ux\}}{I}$, (for a normed field $\k$, complete with respect to its norm);
\item or $R=\quotients{S}{I}$, where $S\sseteq\k[[\ux]]$ is a local regular Henselian, closed under the action of the ordinary partial derivatives, $\di_i(S)\sseteq S$, and $A\in\Mat$ is a matrix of polynomials for some choice of generators;
\item or $R=\quotients{C^\infty(\R^p,0)}{I}$, where the ideal $I$ is finitely generated by some algebraic power  series.
\eee

We remark that among the admissible rings in iii. is e.g. the ring of algebraic power series, $\k<\ux>$, and its quotients. But the localization of the polynomial ring, $\k[\ux]_{(\ux)}$, is not Henselian, thus is not permitted.

\subsubsection{}\label{Sec.Results.Ideals}
We often use the quotient of ideals, $I:J=\{f\in R|\ fJ\sseteq I\}$.

The saturation of $I\sset R$ by $J\sset R$ is the ideal $Sat_J(I):=\suml_{k=1}^\infty I:J^k$, \cite[pg. 318]{Eisenbud}. As $I:J^k\sseteq I:J^{k+1}$ the sum can be thought of as a growing sequence of ideals.
(We do not use the notation $I:J^\infty$ to avoid any confusion with the ideal $J^\infty$ of $J$-flat functions for non-Noetherian rings.)
The saturation can be expressed via the zeroth local cohomology, $Sat_J(I)= H^0_{(J)}(\quotients{R}{I})$, \cite[page 100]{Eisenbud}. Geometrically one erases the subscheme defined by $J$ and then takes the Zariski closure, i.e. $V(Sat_J(I))=\overline{V(I)\smin V(J)}$.

\

Suppose $J\supseteq\cm^\infty$. The Loewy length, $ll_R(J)\le\infty$,  is the minimal number $N\le\infty$
 for which holds $J\supseteq\cm^N$. This number also equals the degree of the socle of the quotient module $\quotients{R}{J}$. For $R=\k[[\ux]]$ we have yet another expression, via the Castelnuovo-Mumford regularity,
 $ll_R(J)=reg(\quotients{R}{J})+1$, see \cite[exercise 20.18]{Eisenbud}.

\

\subsubsection{}\label{Sec.Results.Ann.Coker}
We denote the zero matrix by $\zero$, the identity matrix by $\one$.

For a matrix $A\in\Mat$ denote by $I_j(A)$ the determinantal ideal generated by all the $j\times j$ minors of $A$. (We put $I_0(A)=R$ and $I_{m+1}(A)=\{0\}$.)
 If $A\in\Matm$ then $height(I_{j+1}(A))\le min\big((m-j)(n-j),dim(R)\big)$ and the equality holds generically, see proposition \ref{Thm.background.Fitt.Ideals.generic.height}.

Denote by $ann.coker(A)$ the annihilator-of-cokernel ideal of the homomorphism $R^{\oplus n}\stackrel{A}{\to}R^{\oplus m}$. For the properties/relation/computability of these ideals see \S\ref{Sec.Background.Fitting.Ideals.ann.coker}.

The ideal $ann.coker(A)$ is a refined (partially reduced) version of $I_m(A)$, equivalently: for a module $M$ the ideal $ann(M)$ is a refinement of the minimal Fitting ideal $Fitt_0(M)$. We also use the generalizations,
   $\{ann.coker_j(A)=ann_{j}(coker(A))$, they refine the ideals $I_j(A)$ and $Fitt_{m-j}(coker(A))$, \cite[exercise 20.9]{Eisenbud}. In \S\ref{Sec.Background.Generaliz.ann.coker} we recall the definitions and some relevant  properties.

\

Sometimes we take the $\cm$-adic completion, $R\to\hR$, we denote by $\hA$ the corresponding completion of $A$.

\subsubsection{}

Denote by $Der_\k(R)$ the $R$-module of all the $\k$-linear derivations of the ring. For an ideal $J\sset R$
 denote by $Der_\k(R,J)$ the module of those derivations that send $R$ to $J$. Note that
   $Der_\k(R,J)\supseteq J\cdot Der_\k(R)$, often the two modules have the same rank, $rank(Der_\k(R))$.

For the classical regular rings, e.g. $\k[\![\ux]\!]$, $\k\{\ux\}$, the module $Der_\k(R)$ is generated by the partial derivatives, $Der_\k(R)=\Span_R(\di_1,\dots,\di_p)$.
 Furthermore, in this case $Der_\k(R,J)=J\cdot Der_\k(R)$ and  $rank(Der_\k(R))=dim(R)$.

\

The derivations act on the matrix entry-wise, for any $\cD\in Der_\k(R)$ one has $\cD(A)\in\Mat$. By applying the whole module $Der_\k(R)$ we get the submodule $Der_\k(R)(A)\sseteq \Mat$
and similarly $Der_\k(R,J)(A)\sseteq \Mat$.

Sometimes we need the explicit generating matrix of this submodule, we call it ``the Jacobian matrix of $A$"
\beq
Jac(A):=\{\cD_\al a_{ij}\}_{\substack{(i,j)\in [1..m]\times[1..n]\\\cD_\al\in Der_\k(R)}}.
 \eeq
 Here $\{\cD_\al\}$ are some generators  of $Der_\k(R)$. The matrix $Jac(A)$
  has $mn$ rows (we identify $\Mat= R^{mn}$), while the number of columns depends on $Der_\k(R)$ (and could be infinite).
    This matrix defines the critical locus of the map $Spec(R)\stackrel{A}{\to}Mat_{m\times n}(\k)$, see \cite[Chapter 4]{Looijenga}. While the matrix $Jac(A)$ depends on the choice of generators of $Der_\k(R)$, we will need only the invariants of its image, the module $Der_\k(R)(A)$, these are well defined.

The generating matrix of the module $Der_\k(R,J)(A)$ is the ``$J$-Jacobian matrix":
$Jac^{(J)}(A):=\{\cD_\al a_{ij}\}_{\substack{(i,j)\in [1..m]\times[1..n]\\\cD_\al\in Der_\k(R,J)}}$.

\subsubsection{The singular locus of an ideal}\label{Sec.Results.Sing(J)}
Suppose the expected height of an ideal $J\sseteq\cm$ is  $r$. (The typical example for $J$ is a determinantal ideal.)
We define the singular locus/ideal of $J$ as follows. Fix any system of generators $f_1,\dots,f_N$ of $J$, write them as a column. Applying the derivations of $R$ we get the submodule $Der_\k(R)(\{f_i\})\sseteq R^{\oplus N}$. Then the singular ideal of $J$ is
\beq
Sing(J):=Sing_r(J):=ann_r\quotient{R^{\oplus N}}{J\cdot R^{\oplus N}+Der_\k(R)(\{f_i\})}.
\eeq
Here $ann_r(\dots)$ is the refined version of $Fitt_{N-r}(\dots)$, see \S\ref{Sec.Results.Ann.Coker}.
 Equivalently, the singular ideal is the generalized annihilator-of-cokernel of the large matrix:
\beq\label{Eq.Definition.of.Sing(J)}
Sing_r(J):=ann.coker_r\bpm f_1\dots f_N&\zero&\zero&\{\cD_\al f_1\}\\\zero&f_1\dots f_N&\zero&\{\cD_\al f_2\}\\\dots&\dots&\dots\\
\zero&\zero&f_1\dots f_N&\{\cD_\al f_N\}\epm
\eeq
The ideal $Sing_r(J)$ does not depend on all the choices made, see \S\ref{Sec.Background.Sing(J)} for this and other properties.

Sometimes we need the $\cm$-singular locus, $Sing^{(\cm)}_r(J)$, with $Der_\k(R,\cm)$ instead of $Der_\k(R)$.

In the classical Singularity Theory  the critical/singular loci of a map are often given the Fitting scheme structure, using $I_r(\dots)$ instead of $ann.coker_r(\dots)$. But when working with $T^1_{(\Si,G,A)}$ we need the annihilator-of-cokernel scheme structure.

\

The typical use of this singularity ideal is for the determinantal ideals. For the ideal $J=I_{j+1}(A)$ the expected height is $(m-j)(n-j)$. Fix a column of generators of $I_{j+1}(A)$, say $(\De_1,\dots,\De_N)^t\in R^{\oplus N}$, and apply the derivations to get the module $Der_\k(R)(\{\De_i\})\sseteq R^{\oplus N}$. Then
 \beq
 Sing_{(m-j)(n-j)}(I_{j+1}(A))=ann_{(m-j)(n-j)}\quotient{R^{\oplus N}}{I_{j+1}(A)R^{\oplus N}+Der_\k(R)(\{\De_i\})}.
 \eeq

Often we need only the radical, $\sqrt{Sing(J)}$, then the definition simplifies, in particular instead of the annihilator-of-cokernel one can write the ideal of minors,
$\sqrt{Sing_r(J)}=\sqrt{J+I_r(Der_\k(R)(\{f_i\}))}$.
We remark that $Sing(J)\supseteq J+ann_r\quotients{R^{\oplus r}}{Der_\k(R)(\{f_i\})}$, but the inclusion can be proper.
See \S\ref{Sec.Background.Sing(J)} for all the proofs.

\subsubsection{Generic finite determinacy}
\bed
We say that the generic determinacy holds for a given action $G\circlearrowright\Si$ if for any $A\in\Si$, any number $N<\infty$ and the generic matrix $B\in Mat_{m\times n}(\cm^N)$, such that $A+B\in\Si$, the matrix $A+B$ is finitely determined.
\eed
The generic determinacy implies that the subset of $\Si$ corresponding to the not-finitely-$(\Si,G)$-determined matrices is of infinite codimension in the sense of Mather/Tougeron, see \cite[\S5]{Wall-1981}.

\subsection{Criteria for the $Aut_\k(R)$-action}\label{Sec.Results.Aut(R)}
The group-action $Aut_\k(R)\circlearrowright\Mat$  does not use any matrix structure, thus we can put $m=1$.
Such one-row matrices can be considered as maps, $Spec(R)\stackrel{A}{\to}\k^n$, and in the ``geometric case" (when $Aut_\k(R)=\cR$) we get the classical right equivalence,
$\cR\circlearrowright Maps\Big(Spec(R),(\k^n,0)\Big)$. Therefore the following theorem is the natural extension of numerous classical results.

\bthe\label{Thm.Results.T1.for.Aut(R)} Let $A\in \Si=Mat_{1\times n}(R)$.
\\1.  $ann\big(T^1_{(\Si,Aut_\k(R),A)}\big)=ann.coker(Jac^{(\cm)}(A))$.
\\2. Suppose $R$ has the relevant approximation property (in the sense of \S\ref{Sec.Results.Relevant.Approx})
  and denote $J:=ann.coker(Jac^{(\cm)}(A))\sseteq\cm$. Then for any $B\in Mat_{m\times n}(J\sqrt{J})$ holds: $A+B\stackrel{Aut^{(\sqrt{J})}_\k(R)}{\sim}A$.
  In other words: $Aut_\k(R)(A)\supseteq \{A\}+Mat_{m\times n}(J\sqrt{J})$.
\\3. Suppose $ann.coker(Jac^{(\cm)}(A))\supseteq\cm^\infty$, then the $Aut_\k(R)$-determinacy order  is bounded:
\[
ll_R\Big(ann.coker(Jac^{(\cm)}(A))\Big)-1\le ord^\Si_{Aut_\k(R) }(A)\le ll_R\Big(ann.coker(Jac^{(\cm^2)}(A))\Big)-1.\]
\ethe
(The subgroup $Aut^{(J)}_\k(R)\sseteq Aut_\k(R)$ is defined in example \ref{Ex.Intro.Typical.Groups}.)

\bex\label{Ex.Results.Fin.Determ.Aut(R)}
For $n=1$ the matrix has just one entry, we replace $A$ by $f$. In this case we study the deformations/determinacy of germs of functions.  The matrix $Jac^{(\cm)}(f)$ has one row, thus $ann.coker(Jac^{(\cm)}(f))=I_1(Jac^{(\cm)}(f))$, i.e. the ideal is generated by the entries of this row. Suppose $R$ is regular, then
 $I_1(Jac^{(\cm)}(f))=\cm\cdot(\di_1 f,\dots,\di_p f)$ and Part 2' of theorem \ref{Thm.Results.T1.for.Aut(R)}  gives:
\beq
ll_R\Big(\cm\cdot(\di_1 f,\dots,\di_p f)\Big)-1\le ord_{Aut_\k(R) }(f)\le ll_R\Big(\cm^2\cdot(\di_1 f,\dots,\di_p f)\Big)-1.
\eeq
For $R=\k[\![\ux]\!]$ or $\k\{\ux\}$, when these bounds are finite, we get e.g. \cite[Part 1 of Theorem 2.23]{GLS}, \cite[Theorem 1.2]{Wall-1981}.
If $R=C^\infty(\R^p,0)$ and $I_1(Jac^{(\cm)}(f))\supseteq\cm^\infty$ this gives the infinite determinacy of \cite[Theorem 6.1]{Wall-1981}. If the ideal $I_1(Jac^{(\cm)}(f))$ contains no $\cm^N$ (for any $N<\infty$) then there is no finite determinacy. See example \ref{Ex.Proofs.Relative.Determinacy.Aut(R)} for the `admissible' deformations and the `relative' determinacy.
\eex

In the case $n>1$ the finite $Aut_\k(R)$-determinacy of maps is a very restrictive condition, it means  the simultaneous $Aut_\k(R)$-determinacy of a tuple of functions.
\bprop\label{Thm.Finite.Determinacy.Coord.Change} Let $R$ be a local ring with the relevant approximation (in the sense of \S\ref{Sec.Results.Relevant.Approx}), and suppose $Der_\k(R)$ is a free module of rank=$dim(R)$.
Suppose $n>1$.
 Then the map $A\in Maps(Spec(R),(\k^n,0))$ is finitely-$Aut_\k(R)$-determined iff
 $dim(R)\ge n$ and the entries of $A$ form a sequence of generators of $\cm$ (over $R$).
 In the later case $A$ is $Aut_\k(R)$-stable.
\eprop
(The condition on $Der_\k(R)$ is satisfied e.g. for any regular local subring of $\k[[\ux]]$.)
 For $R=\C\{\ux\}$ this was proved by Mather, see e.g. proposition 2.3 in \cite{Wall-1981}.

The proofs and examples are in \S\ref{Sec.Proofs.Aut(R)}.

\subsection{Criteria for the $\cG_{lr}$ action}\label{Sec.Results.Criteria.G.coord.with.changes}
First we study the annihilator of the module $T^1_{(\Mat,\cG_{lr},A)}$.
\bthe\label{Thm.Results.annT1.Grl.coord.changes}
Let $R$ be a local Noetherian ring and $A\in \Si=\Mat$.
\\1. If $m=1$ then $ann\big(T^1_{(\Si,\cG_{lr},A)}\big)=Sing^{(\cm)}_{n}(I_1(A))\supseteq I_1(A)+ann.coker(Jac^{(\cm)}(A))$.
\\2. For any non-maximal prime ideal, $\cp\ssetneq\cm$, satisfying $\cp\supseteq I_m(A)$ but $\cp\not\supseteq I_{m-1}(A)$, there holds:
\[ann(T^1_{(\Si,\cG_{lr},A)})_{\cp}=Sing_{n-m+1}(I_m(A))_{\cp}.\]
\\3.
$ann.coker(A)+ann.coker(Jac^{(\cm)}(A))\sseteq ann(T^1_{(\Si,\cG_{r},A)})\sseteq ann(T^1_{(\Si,\cG_{lr},A)})\sseteq
\capl^{m-1}_{j=0}Sat_{I_j(A)}\Big(Sing^{(\cm)}_{(m-j)(n-j)}\big(I_{j+1}(A)\big)\Big)$.
\\3'. In particular, if $rank_R(Der_\k(R))<(m-j)(n-j)$ or $height(I_{j+1}(A))<(m-j)(n-j)$, for some $j$, then $ann(T^1_{(\cG_{lr},\Si,A)})\sseteq Sat_{I_j(A)}(I_{j+1}(A))$.
\\4. If $R$ is Noetherian then $\sqrt{ann(T^1_{(\Si,\cG_{lr},A)})}=
\sqrt{\capl^{m-1}_{j=0}Sat_{I_{j}(A)}\Big(Sing_{(m-j)(n-j)}\big(I_{j+1}(A)\big)\Big)}$.
\ethe
In part 1. the matrix structure plays no role and (when $Aut_\k(R)=\cR$) the $\cG_{lr}$-action induces the classical contact equivalence of maps, $\cK$. Accordingly part 1 extends the classical criteria of determinacy of maps, see  \S\ref{Sec.Results.Contact.Equiv}.

\

We give the geometric interpretation of the other statements. Note that the support of $T^1_{(\Si,\cG_{lr},A)}$ lies inside the locus $V(I_m(A))$. (In fact, if $A$ is non-degenerate at some point then it is $G_{lr}$-stable at that point.) Part 2 says that the support of $T^1_{(\Si,\cG_{lr},A)}$ contains the singular locus of $V(I_m(A))$, with the locus $V(I_{m-1}(A))$ erased. The embedding $ann(T^1_{(\Si,\cG_{lr},A)})\sseteq\dots$ of part 3 reads:
\beq
supp\Big(T^1_{(\Si,\cG_{lr},A)}\Big)\supseteq \cupl^{m-1}_{j=0} \overline{Sing\Big(V\big(I_{j+1}(A)\big)\Big)\smin V\big(I_{j}(A)\big)}.
\eeq
(This is the embedding of schemes, the closure is taken in Zariski topology.) Finally, part 4 is a set-theoretic equality for the reduced schemes:
\beq
supp\Big(T^1_{(\Si,\cG_{lr},A)}\Big)_{red}=\Big(\cupl^{m-1}_{j=0} \overline{Sing\Big(V\big(I_{j+1}(A)\big)\Big)\smin V\big(I_{j}(A)\big)}\Big)_{red}.
\eeq
More detail and a separate proof of the set-theoretic equality are in \S\ref{Sec.Proof.Supp(T1).Glr.set.theoretic}.

\

\

Suppose $ann(T^1_{(\Si,\cG_{lr},A)})\supseteq\cm^\infty$, then the order of determinacy is bounded by the Loewy length of the annihilator (see proposition \ref{Thm.Background.Ord.of.det.via.Loewy.length}):
\beq
ll_R(ann(T^1_{(\Si,\cG_{lr},A)}))-1\le ord^\Si_{\cG_{lr}}(A)\le
ll_R(ann(T^1_{(\Si,\cG^{(\cm)}_{lr},A)}))-1.
\eeq

We give some immediate applications to the finite determinacy of matrices.
\bprop\label{Thm.Results.Finite.Determinacy.Glr.with.coord.change.Corol} Let $R$ be a local Noetherian ring.
\\1. Suppose  $rank(Der_\k(R))<(n-j)(m-j)\le dim(R)$ for some $1\le j< m$. Then
 no $A\in\Matm$ is finitely-$\cG_{lr}$-determined.
\\2. If $A\in\Matm$ is finitely-$\cG_{lr}$-determined then all the ideals $\{I_j(A)\}$ are of expected heights.
\\For the rest of the proposition suppose that $R$ has the relevant approximation property (in the sense of \S\ref{Sec.Results.Relevant.Approx}).
\\3. If $dim(R)\le n-m+1$ then $A\in\Matm$ is finitely-$\cG_{lr}$-determined iff it is finitely-$G_r$-determined.
\\4. Fix $A\in \Mat$ and let $J:=ann(T^1_{(\Si,\cG_{lr},A)})\sseteq\cm$.
Then for any $B\in Mat_{m\times n}(J\cdot\sqrt{J})$ holds:
 $A+B\stackrel{\cG^{(\sqrt{J})}_{r}}{\sim}A$.
\eprop

\bex
\bee
\item Part 4 of the proposition strengthens and extends Theorem 2.1 of \cite{Cutkosky-Srinivasan}, which was proven for the case  $m=1$, $R=\k\{\ux\}$ and with $I_{n}(Jac^{(\cm)}(A))$ instead of $ann.coker(Jac^{(\cm)}(A))$.

\item Geometrically part 2 of the proposition reads, for $\k=\bar\k$: if $A$ is finitely-$\cG_{lr}$-determined then all the loci $V(I_j(A))\sset Spec(R)$ are of the expected codimensions.
 In particular:
\li either $dim(R)>mn$ and $I_1(A)$ is a complete intersection ideal
\li or $dim(R)\le mn$ and $I_1(A)$ contains a power of the maximal ideal.

\item Suppose $\k=\bar\k$ and the ring $R$ is regular (i.e. $Spec(R)$ is smooth), then $A$ is finitely determined iff
all the loci $V(I_j(A))\sset Spec(R)$ are of the expected codimensions and
all the complements $V(I_j(A))\smin V(I_{j-1}(A))$ are smooth. For square matrices over $\C\{\ux\}$ this is proposition 3.2 of \cite{Bruce-Tari04}.

\item  The theorem implies: for regular rings (with relevant approximation properties)
 the generic $\cG_{lr}$-finite-determinacy holds. This extends the classical Tougeron/Mather criteria, see \S\ref{Sec.Results.Relation.to.Singular.Theory}. (Again, for square matrices over $\C\{\ux\}$ this is proposition 3.3 of \cite{Bruce-Tari04}.)
\eee
\eex

\

\bprop\label{Thm.Intro.Glr.coord.changes.Additional.Corollary}
  Suppose $\k=\bar\k$, $R$ is local Noetherian, regular, with the relevant approximation property (in the sense of \S\ref{Sec.Results.Relevant.Approx}).
\\1. If $dim(R)\le 2(n-m+2)$, then $A$ is finitely-$\cG_{lr}$-determined iff
the locus $V(I_{m}(A))\sset Spec(R)$ is of expected codimension and $Sing_{n-m+1}(V(I_m A))=V(I_{m-1}(A))=\{0\}\sset Spec(R)$. (``set-theoretically", i.e. with the reduced scheme structure)
\\2. In particular, if $dim(R)=n-m+2$, then $A$ is finitely-$\cG_{lr}$-determined iff $I_m(A)$
defines an isolated reduced curve singularity.
\\3. If $A$ is finitely-$\cG_{lr}$-determined then $I_1(A)$ defines a subgerm of $Spec(R)$ with at most an isolated singularity.
\eprop

The proofs, examples and further applications to the finite determinacy are in \S\ref{Sec.Proofs.cG_lr}

\subsection{Applications to the finite determinacy of maps/complete intersections under the contact equivalence}\label{Sec.Results.Contact.Equiv}
In this section we assume that $R$ has the approximation property in the sense of \S\ref{Sec.Results.Relevant.Approx} and moreover $Aut_\k(R)=\cR$. (The later holds e.g. for any subring of $\quotients{\k[[\ux]]}{I}$, see  \cite[\S3.1]{Belitski-Kerner}.)

For matrices with just one row, $Mat_{1\times n}(R)$, the $\cG_{lr}$-equivalence coincides with the contact equivalence of maps/complete intersections, $\cK\circlearrowright Maps(Spec(R),(\k^n,0))$, see \S\ref{Sec.Results.Relation.to.Singular.Theory}. Therefore we replace $A$ by $f$.  In this case part 1 of theorem
 \ref{Thm.Results.annT1.Grl.coord.changes} gives the exact expression for the annihilator and the
 criteria are especially simple.
\bex
In the simplest case, $m=1=n$, we get the contact-determinacy of the function germs/hypersurface singularities.
 Here  $Sing^{(\cm)}(I_1(f))=I_1(f)+ann.coker(Jac^{(\cm)}(f))$.
 Suppose $R$ is regular, then $Sing^{(\cm)}(I_1(f))=(f)+\cm(\di_1 f,\dots,\di_p f)$.
\begin{enumerate}[i.]
\item If this ideal contains a power of the maximal ideal then we get Part 2 of \cite[Theorem 2.23]{GLS}:
  \beq
  ll_R\Big((f)+\cm(\di_1 f,\dots,\di_p f)\Big)-1\le ord_{\cK}(f)\le ll_R\Big(\cm(f)+\cm^2(\di_1 f,\dots,\di_p f)\Big)-1.
  \eeq
   If the ideal $(f)+\cm(\di_1 f,\dots,\di_p f)$ contains (at least) $\cm^\infty$ then we get the infinite determinacy, \cite[theorem 6.1]{Wall-1981}.
\item If $R$ is regular then by the Brian\c{c}on-Skoda theorem $f^{dim(R)}\in Jac(f)$. (The initial version was for analytic rings, for the general result see \cite{Lipman-Sathaye}.) In particular, $f$ is finitely-$\cK$-determined iff $f$ is finitely-$\cR$-determined. For $R=\C\{\ux\}$ this is proposition 2.3 of \cite{Wall-1981}.

    More generally, this relation of finite determinacies holds for any $R$ and $f$ satisfying $f\in \sqrt{Jac(f)}$.
\item  When the ideal $Sing(I_1(f))$ contains no power of $\cm$ there is no finite determinacy. (For $\k=\bar\k$ this corresponds to the non-isolated singularity.) Suppose $R$ is Noetherian and
 define $J=\sqrt{(f)+ann.coker(Jac^{(\cm)}(f))}$. This ideal defines (set-theoretically) the singular locus of the hypersurface $\{f=0\}\sset Spec(R)$. As the singularity is non-isolated $J\subsetneq\cm$.
 We want to understand the ``space of admissible deformations", i.e. the biggest ideal $I\sset R$ such that $f$ is finitely determined for deformations inside $\Si=I$. By proposition \ref{Thm.Background.Reduction.to.annihilator} we have
 \beq
J\cdot Sat_\cm\Big((f)+ann.coker(Jac^{(\cm)}(f))\Big)\sseteq  I\sseteq Sat_\cm \Big((f)+ann.coker(Jac^{(\cm)}(f))\Big)
 \eeq
 \end{enumerate}
\eex
\bex
Now consider the case $m=1<n$.
\bee[i.]
\item
We get: $ll_R(Sing^{(\cm)}(\uf))-1\le ord_\cK(\uf)\le ll_R(Sing^{(\cm^2)}(\uf))-1$. This generalizes e.g. proposition 4.1 of \cite{Wall-1981}.
\item
If $dim(R)\le n$ then $\uf$ is finitely-$\cK$-determined iff it is finitely-$G_r$-determined.
 For $dim(R)\le n$  the generic matrix is finitely determined.
 \item
 $\uf$ is finitely-$\cK$-determined iff  $Sing(\uf)\supseteq\cm^N$, for some $N$. $\uf$ is infinitely-$\cK$-determined iff  $Sing(\uf)\supseteq\cm^\infty$.
 \item
If $R$ is regular and $dim(R)>n$ then $\uf$ is finitely-$\cK$-determined iff
 $\uf^{-1}(0)$ is a complete intersection (of codimension $n$) with at most an isolated singularity.
 In particular, the generic map $Spec(R)\stackrel{\uf}{\to}(\k^n,0)$ is finitely-$\cK$-determined.
\eee
\eex

\beR\label{Ex.Fin.Det.Nonisolated.Sing.for.non-regular.ring}
Part (iv.) of the corollary is the classical criterion, e.g. \cite{Wall-1981}.  For non-regular rings the locus
 $\uf^{-1}(0)$ can have non-isolated singularity but $\uf$ can still be finitely determined. For example,
 let $R=\k[\![x_1,\dots,x_p]\!]/(x^k_1)$, then the module of derivations is generated by $(x_1\di_1,\di_2,\dots,\di_p)$.
 Consider $\uf=x_1+x_p\in Mat_{1\times 1}(\cm)$. Then $T_{(\cG_r\uf,\uf)}=T_{(Aut_\k(R)\uf,\uf)}=\cm$, thus $ann(T^1_{(R,Aut_\k(R),\uf)})=\cm$,
 implying the finite-$Aut_\k(R)$-determinacy, hence the finite-$\cK$-determinacy.
  But the zero locus is
 $\{x_1+x_n=0\}\approx Spec(\k[x_1,x_2,\dots,x_{n-1}]/(x^k_1))$, i.e. is a multiple hyperplane.
\eeR

\subsection{Congruence and (skew-)symmetric matrices}\label{Sec.Results.Congruence}
 In \cite{Belitski-Kerner1} we have studied the actions $G_{congr}\circlearrowright \Si=Mat_{m\times m}(R)$, $G_{congr}\circlearrowright \Si^{sym}=Mat^{sym}_{m\times m}(R)$, $G_{congr}\circlearrowright \Si^{skew-sym}=Mat^{skew-sym}_{m\times m}(R)$. We have shown (in parts 2,3 we assume that $R$ is Henselian and Noetherian):
\bee
\item   $ann\Big( T^1_{(\Si,G_{congr},A)}\Big)\supseteq\cm^\infty$. In particular, the module $T^1_{(\Si,G_{congr},A)}$ is supported generically on $Spec(R)$ and
 there are no finitely-$(\Si,G_{congr})$-determined matrices in $Mat_{m\times m}(R)$.
\item The height of $ann\Big( T^1_{(\Si^{sym},G_{congr},A)}\Big)$ is at most 1. If $dim(R)>1$ then no $A\in Mat^{sym}_{m\times m}(\cm)$ is finitely determined.
If $dim(R)\le1$ then
$ll_R\Big(\overline{I_m(A):I_{m-1}(A)}\Big)-1\le ord^{\Si^{sym}}_{G_{congr}}(A)\le ll_R\Big(I_m(A):I_{m-1}(A)\Big)$.
 In particular the finite-$(\Si^{sym},G_{congr})$-determinacy holds generically.

The same holds also for $ann\Big( T^1_{(\Si^{skew-sym},G_{congr},A)}\Big)$, when $m$ is even.
\item Suppose $m$ is odd,  then the height of $ann\Big( T^1_{(\Si^{skew-sym},G_{congr},A)}\Big)$ is at most 3.  If $dim(R)>3$ then no $A\in Mat^{skew-sym}_{m\times m}(\cm)$ is finitely-$(\Si^{skew-sym},G_{congr})$-determined. If $dim(R)\le3$ then
\[
ll_R\Big(\overline{I_{m-1}(A):I_{m-2}(A)}\Big)-1\le ord^{\Si^{skew-sym}}_{G_{congr}}(A)\le ll_R(Pf_{m-1}(A)).
\]
(Here $Pf_{m-1}(A)$ is the sum of Pfaffian ideals, $\sum Pf_{m-1}(A_{ii})$, of the $(m-1)\times(m-1)$ blocks, obtained by erasing $i$'th row and column, for $1\le i\le m$.)
 In particular, the finite-$(\Si^{skew-sym},G_{congr})$-determinacy holds generically.
\eee

\

The following theorem addresses the larger group, $\cG_{congr}$.
\bthe\label{Thm.Results.T1.for.Congruence}
Let $R$ be a local Noetherian ring, $dim(R)>0$.
\\1. Suppose $\cm^\infty=\{0\}$ and $rank(Der_\k(R))<\lfloor\frac{m}{2}\rfloor$ then $ann\Big( T^1_{(Mat_{m\times m}(R),\cG_{congr},A)}\Big)=\{0\}$.
\\2. Fix some $A\in \Si=Mat^{sym}_{m\times m}(R)$.
\bee[i.]
\item
For any non-maximal prime ideal, $\cp\ssetneq\cm$, satisfying $\cp\supseteq (det(A))$ but $\cp\not\supseteq I_{m-1}(A)$, there holds:
\[
(T^1_{(\Si,\cG_{congr},A)})_{\cp}=\quotient{R_\cp}{\Big((det(A)_\cp+ Jac^{(\cm)}(det(A)_\cp)\Big)},\quad \text{ in particular }
ann(T^1_{(\Si,\cG_{congr},A)})_{\cp}=det(A)_{\cp}+I_1\Big(Jac^{(\cm)}(det(A))\Big).\]
\item
$ann.coker(A)+ann.coker(Jac^{(\cm)}(A))\sseteq  ann(T^1_{(\Si,\cG_{congr},A)})\sseteq
\capl^{m-1}_{j=0}Sat_{I_j(A)}\Big(Sing_{\bin{m-j+1}{2}}(I_{j+1}(A))\Big)$.
\item If  $rank_R(Der_\k(R))<\bin{m-j+1}{2}$ or  $height(I_{j+1}(A))<\bin{m-j+1}{2}\le dim(R)$, for some $0\le j<m$, then
\[
ann(T^1_{(\Si,\cG_{congr},A)})\sseteq Sat_{I_j(A)}(I_{j+1}(A)).
\]
\item  $\sqrt{ann\Big(T^1_{(\Si,\cG_{congr},A)}\Big)}=
\capl^{m-1}_{j=0}\sqrt{Sat_{I_j(A)}\Big(Sing_{\bin{m-j+1}{2}}(I_{j+1}(A))\Big)}$.
\eee
3. Let $A\in \Si=Mat^{skew-sym}_{m\times m}(R)$.
\bee[i.]
\item Suppose $m$ is even and a non-maximal prime ideal, $\cp\ssetneq\cm$, satisfies $\cp\supseteq (det(A))$ but $\cp\not\supseteq I_{m-2}(A)$. Then
$ann(T^1_{(\Si,\cG_{congr},A)})_{\cp}=Sing_{1}(I_{m-1}(A))_{\cp}$.
\item
For any $m$ (even or odd) holds:
$\sqrt{ann\Big(T^1_{(Mat^{skew-sym}_{m\times m}(R),\cG_{congr},A)}\Big)}=
\capl^{m-1}_{\substack{j=0\\j\in2\Z}}\sqrt{Sat_{I_j(A)}\Big( Sing_{\bin{m-j}{2}}\big(I_{j+1}(A)\big)\Big)}$.
\eee
\ethe
Parts 2,3 of the theorem have the immediate set-theoretic interpretations, similar to those of theorem \ref{Thm.Results.annT1.Grl.coord.changes}.

The following proposition shows how the increasing the group improves the determinacy.
\bprop\label{Thm.Results.Congruence.Finite.Determin}
1. If $rank(Der_\k(R))<\lfloor\frac{m}{2}\rfloor$ then no $A\in Mat_{m\times m}(R)$ is finitely $\cG_{congr}$-determined.
\\2. Suppose  $rank(Der_\k(R))<\bin{m-j+1}{2}\le dim(R)$ for some $0\le j< m$. Then
 no $A\in Mat^{sym}_{m\times m}(\cm)$ is  finitely-$\cG_{congr}$-determined.
\\For the rest of the proposition suppose that $R$ has the relevant approximation property (in the sense of \S\ref{Sec.Results.Relevant.Approx}).
\\3.  If $A\in Mat^{sym}_{m\times m}(\cm)$ is finitely-$\cG_{congr}$-determined then all the ideals $\{I_j(A)\}$ are of expected heights.
\\4. If $dim(R)\le1$ then $A\in Mat^{sym}_{m\times m}(\cm)$ is finitely-$\cG_{congr}$-determined iff it is finitely-$G_{congr}$-determined.
\\5. Fix $A\in \Si=Mat^{sym}_{m\times m}(\cm)$ and denote $J:=ann(T^1_{(\Si,\cG_{congr},A)})\sseteq\cm$.
Then for any $B\in Mat_{m\times m}^{sym}(J\cdot\sqrt{J})$ holds:
 $A+B\stackrel{\cG^{(\sqrt{J})}_{congr}}{\sim}A$.
\eprop
Similar statements hold for $A\in Mat^{skew-sym}_{m\times m}$.

\subsection{Remarks}\label{Sec.Intro.Remarks}
\subsubsection{``Complexity of the answers''}
As one sees, the bounds on $ann(T^1_{(\Si,G,A)})$ are somewhat involved.
In the (classical/simplest) case of functions with isolated singularities, $m=1=n$, the tangent spaces to the miniversal deformation are quite simple,
 $T^1_{(Aut_\k(R),f)}=\quotients{R}{Jac_f}$,  $T^1_{(\cK,f)}=\quotients{R}{(f)+Jac_f}$. Already the complete intersection case, $m=1\le n$, is cumbersome.
 Therefore, one cannot hope for a nice/short expression of $ann(T^1_{(\Si,G,A)})$ in the matrix case. The best we can try is to express the upper/lower bounds on $ann(T^1_{(\Si,G,A)})$ in terms of some invariants of $A$ and to interpret these geometrically.

Thus theorems  \ref{Thm.Results.T1.for.Aut(R)}, \ref{Thm.Results.annT1.Grl.coord.changes},
\ref{Thm.Results.T1.for.Congruence} and their applications to the finite determinacy go in the spirit
 of the classical Mather's criterion that transforms  (for $\k=\bar\k$) the general statement

$\bullet$ "the hypersurface germ is finitely determined iff its miniversal deformation is finite dimensional".
\\to the geometric condition

$\bullet$ "the hypersurface germ is finitely determined iff it has at most an isolated singularity".

\

\subsubsection{Relations to Singularity Theory}\label{Sec.Results.Relation.to.Singular.Theory}
For a short discussion of the results on determinacy see \cite[\S2.9]{Belitski-Kerner1}
\li The case of square matrices (for $\k=\R$ or $\k=\C$, $R=\k\{x_1,\dots,x_p\}$ and $G=\cG_{lr}$)
was considered in \cite{Bruce-Tari04}, and further studied in \cite{Bruce-Goryun-Zakal02},
\cite{Bruce2003}, \cite{Goryun-Mond05}, \cite{Goryun-Zakal03}. In particular, the generic finite
determinacy was established and the simple types were classified.
\li Finite determinacy is equivalent to the finite dimensionality of the miniversal deformation.
 In particular, the genericity of finite determinacy, for a fixed action $G\circlearrowright\Si$, means:
 the stratum of matrices with Tjurina number $\infty$, $\{A\}+\Si^G_{\tau=\infty}\sset \{A\}+\Si$, is of infinite codimension.

\subsubsection{Transition from $\cG_{lr}$ to $\cK$ and back}
One often forgets the matrix structure and considers matrices as  maps, $Maps(Spec(R),(\k^{mn},0))$.
 This increases the group, from $\cG_{lr}=GL(m,R)\times GL(n,R)\rtimes Aut_\k(R)$ to $GL(mn,R)\rtimes Aut_\k(R)$. When $Aut_\k(R)=\cR$ the equivalence induced by  $GL(mn,R)\rtimes Aut_\k(R)$ coincides with the contact equivalence of maps, $\cK$.

In the converse direction, we can start from the space $Maps(Spec(R),(\k^N,0))$, with $N=mn$, $1<m\le n$, and identify $\k^{mn}\approx Mat_{m\times n}(\k)$.
 This associates to any map the corresponding matrix $A^{(m,n)}\in\Matm$.
This decreases the group, from $GL(mn,R)\rtimes Aut_\k(R)\circlearrowright Maps\big(Spec(R),(\k^N,0)\big)$ to the subgroup  $\cG^{(m,n)}_{lr}:=GL(m;R)\times GL(n;R)\rtimes Aut_\k(R)\circlearrowright\Matm$.  The orbits of $\cG^{(m,n)}_{lr}$  are {\em much smaller} than those of $GL(mn,R)\rtimes Aut_\k(R)$. (The $\cG^{(m,n)}_{lr}$-orbits are of infinite codimension inside  the $GL(mn,R)\rtimes Aut_\k(R)$-orbits.)
 Thus, the finite $\cG^{(m,n)}_{lr}$-determinacy is a much stronger
  property than the finite $GL(mn,R)\rtimes Aut_\k(R)$-determinacy.
For $m=n$ we can decrease the group further, to $\cG_{congr}$.

 From our results we see: even for the group $\cG^{(m,n)}_{lr}$ the finite determinacy is often the generic property. (At least this holds when $R$ is a regular ring). This both strengthens Mather's/Tougeron's results and extends them to the broader category.

\section{Preparations}\label{Sec.Preparations}

\subsection{The relevant algebra-geometric dictionary}
While in this paper we often work in the generality ``$R$ is a local Noetherian ring", our statements and proofs usually have very transparent geometric/set-theoretic interpretation. (In fact the proofs were initially obtained ``set-theoretically" and then have been translated to the algebraic language.) In this section  we collect several results of algebraic-geometric dictionary. Though these are standard results in commutative algebra, we give the proofs and the geometric/set-theoretic interpretations, replacing the ring $R$ by the germ of the space $Spec(R)$, the ideal $J$ by the corresponding locus $V(J)$ etc.

For the set-theoretic interpretation we usually assume: $R=\C\{x_1,\dots,x_p\}$ and $Spec(R)$ denotes the germ $(\C^p,0)$ with continuum of reduced/closed points near the origin.

\subsubsection{Points in the neighborhood of the origin vs localizations at prime ideals}\label{Sec.Geometric.Rings.of.Functions}
Sometimes $R$ is the ring of ``genuine" functions, i.e. for any element $f\in R$ the germ $Spec(R)$ has a
representative that contains other closed points besides the origin and $f$ can be actually computed at those
points ``off the origin".
Thus, for any $A\in\Mat$ we can take a small enough representative of $Spec(R)$
 and for any point of it we can evaluate the matrix, $A|_{pt}\in Mat_{m\times n}(\k)$.

For example, this happens for the ring of rational functions on $(\k^p,0)$ that are regular at the origin. More generally, this holds for  the localization at the origin of an affine ring. Another example is the ring of converging power series, $\C\{\ux\}$.
Geometry always suggests how to formulate criteria. Usually the
geometric conditions are of the type: a property $\mathcal{P}$ is satisfied ``generically"
 on  $Spec(R)$; or is satisfied on $Spec(R)\setminus\{0\}$; or on some locus in $Spec(R)$.

 When $R$ is the ring of ``genuine" functions this means:
 for a small enough representative $\cU$ of $Spec(S)$, there exists a locus $\cV\sset\cU$ (open or locally closed)
  such that the condition $\mathcal{P}$  is satisfied at each point of $\cV$.

Complete rings, e.g. $\k[\![\ux]\!]$, are not ``geometric", their elements in general cannot be computed ``off the origin". In this case the language of points is replaced by the language of localizations at prime ideals.

Finally, we emphasize that when speaking (geometrically/set-theoretically) about the points near the origin we assume the field is algebraically closed, $\k=\bar\k$.

\subsubsection{Some generalities on the localizations}
\bel\label{Thm.Background.Localization.1}
Let $R$ be a local Noetherian ring and $\cp\sseteq\cm$ a prime ideal.
\\0. $\cp\not\supseteq J$ iff $J_\cp=R_\cp$.
\\1. If $\cp\not\supseteq I$ then $Sat_I(J)_\cp=J_\cp$. (The saturation is defined in \S\ref{Sec.Results.Ideals}.)
\\2. Fix some ideals $I,J_1,J_2\sset R$, then
 $Sat_I(J_1)\sseteq Sat_I(J_2)$ iff for any prime $\cp\not\supseteq I$ holds: $(J_1)_\cp\sseteq(J_2)_\cp$.
\eel
\bpr
0.  $\Rrightarrow$ If $\cp\not\supseteq J$ then exists an element $f\in J\smin\cp$. As $\cp$ is prime, $f$ is not a zero divisor on $R\smin\cp$. Thus the image of $f$ in the localization $J_\cp$ is invertible. Hence $J_\cp=R_\cp$.

$\Lleftarrow$ If $\cp\supseteq J$ then $J_\cp\sseteq \cp_\cp\ssetneq R_\cp$.

1.  $\supseteq$ is obvious as $Sat_I(J)\supseteq J$.  For the part $\sseteq$ suppose $f\in Sat_I(J)$ then $I^Nf\sseteq J$ for some $N$. As $\cp\not\supseteq I$ there exists an element $g\in I^N$ such that $g\not\in\cp$. But then the image of $g$ in $R_\cp$ is invertible. Hence $(f)_\cp\sseteq J_\cp$.

2. $\Rrightarrow$ if $\cp\not\supseteq I$ then by part (1): $Sat_I(J_i)_\cp=(J_i)_\cp$, thus $(J_1)_\cp\sseteq (J_2)_\cp$.

$\Lleftarrow$ Consider the quotient $\quotients{Sat_I(J_1)+Sat_I(J_2)}{Sat_I(J_2)}$. The localization of this quotient at any prime $\cp\not\supseteq I$ vanishes:
 \beq
\Big(\quotient{Sat_I(J_1)+Sat_I(J_2)}{Sat_I(J_2)}\Big)_\cp=
\quotient{Sat_I(J_1)_\cp+Sat_I(J_2)_\cp}{Sat_I(J_2)_\cp}\stackrel{Part\ 1}{=}
\quotient{(J_1)_\cp+(J_2)_\cp}{(J_2)_\cp}=\{0\}.
 \eeq
 Therefore this quotient is not supported at any $V(\cp)$ for $\cp\not\supseteq I$. Thus there exists $N$ satisfying:
 \beq
 ann(\quotient{Sat_I(J_1)+Sat_I(J_2)}{Sat_I(J_2)})\supseteq I^N\quad\rightsquigarrow\quad
 I^N\cdot Sat_I(J_1)\sseteq Sat_I(J_2)\quad\rightsquigarrow\quad
 Sat_I(J_1)\sseteq Sat_I(J_2).\bull
 \eeq

{\em \underline{The geometric interpretation:}}
\\0. $x\not\in V(J)$ iff the germ $(V(J),x)$ is empty.
\\1. Suppose $x\not\in V(I)$. Then $x\in\overline{V(J)\smin V(I)}$ iff $x\in V(J)$.
\\2. $\overline{V(J_1)\smin V(I)}\supseteq \overline{V(J_2)\smin V(I)}$ iff for any point $x\not\in V(I)$ holds:
 $(V(J_1),x)\supseteq (V(J_2),x)$.

\

\bel\label{Thm.Background.Localization.2}
Let $R$ be a local Noetherian ring and  $J_1,J_2\sseteq \cm$ some  proper ideals. The following conditions are equivalent:
\\1. $\sqrt{J_1}=\sqrt{J_2}$.
\\2. For any non-maximal prime ideal, $\cp\ssetneq\cm$, holds: $(J_1)_\cp\neq R_\cp$ iff $(J_2)_\cp\neq R_\cp$.
\\3.  For any non-maximal prime ideal, $\cp\ssetneq\cm$, holds: $(J_1)_\cp= R_\cp$ iff $(J_2)_\cp= R_\cp$.
\eel
\bpr
Obviously $2\Leftrightarrow3$, thus we prove $1\Leftrightarrow3$.

$1\Rightarrow3$ If $(J_1)_\cp=R_\cp$ then there exists $f\in J_1$ whose image in $(J_1)_\cp$ is invertible. Thus $f\not\in\cp$. But $f^N\in J_2$ for some $N<\infty$. And the image of $f^N$ in $(J_2)_\cp$ is still invertible, hence $(J_2)_\cp=R_\cp$.

$3\Rightarrow1$ Take the prime decomposition, $\sqrt{J_1}=\cap \cp_i$.
(As $\sqrt{J_1}$ is a radical ideal, its primary decomposition consists of prime ideals.)
Suppose for some $i$ happens $\cp_i\not\supseteq\sqrt{J_2}$ then $\cp_i\not\supseteq J_2$, thus $(J_2)_{\cp_i}=R_{\cp_i}\supsetneq (J_1)_{\cp_i}$. Thus, if $\cp$ is a minimal prime for $J_1$ then $\cp_i\supseteq J_2$. Suppose $\cp_i$ is not a minimal prime for $J_2$, then exists a smaller prime ideal $\cq\ssetneq\cp_i$, which is a minimal prime for $J_2$. But then, by the same argument as above, $\cq\supseteq J_1$, thus $\cp_i$ could not be a minimal prime.

Therefore: $\cp$ is a minimal prime ideal for $\sqrt{J_1}$ iff it is the one for $\sqrt{J_2}$.
In other words, $\sqrt{J_1},\sqrt{J_2}$ have the same minimal primes.
\epr

{\em \underline{The geometric interpretation:}} Denote by $0\in Spec(R)$ the base point of the germ.
 Two proper ideals define (set-theoretically) the same locus, $V(J_1)_{red}=V(J_2)_{red}$, iff for any point, $0\neq pt\in Spec(R)$, there holds: $pt\in V(J_1)$ iff $pt\in V(J_2)$.

\subsection{Determinantal ideals and annihilator of cokernel}\label{Sec.Background.Fitting.Ideals.ann.coker}
\cite[\S 20]{Eisenbud}
For $1\le j\le m$ and $A\in Mat_{m\times n}(R)$ denote by $I_j(A)\sset R$
the determinantal ideal generated by all the $j\times j$ minors of $A$.
By definition $I_0(A)=R$ and $I_{>m}(A)=\{0\}$. When the size of $A$ is not stated explicitly we denote the ideal of maximal minors by $I_{max}(A)$.
 The chain of ideals  $R=I_0(A)\supseteq I_1(A)\supseteq\cdots\supseteq I_m(A)$ is {\em invariant} under the $G_{lr}$-action. Moreover, any element $\phi\in Aut_\k(R)$ induces the map of chains,
\beq
\Big\{\cdots\sseteq I_j(A)\sseteq I_{j-1}(A)\sseteq\cdots\Big\}\stackrel{\phi}{\to}
\Big\{\cdots\sseteq I_j(\phi(A))\sseteq I_{j-1}(\phi(A))\sseteq\cdots\Big\}
\eeq
The height of $I_j(A)$ is at most $min\Big((m+1-j)(n+1-j),dim(R)\Big)$ and the equality holds
 generically.
 More precisely:
\bprop\label{Thm.background.Fitt.Ideals.generic.height}
For any $A\in Mat_{m\times n}(\cm)$, for any $N>0$ and the generic $B\in Mat_{m\times n}(\cm^N)$,
  the height of $I_j(A+B)$ is $min\Big((m+1-j)(n+1-j),dim(R)\Big)$.
\eprop

\

A matrix $A\in\Matm$ can be considered as a homomorphism of free $R$-modules, its cokernel is an $R$-module as well:
\beq
 R^{\oplus n}\stackrel{A}{\to} R^{\oplus m}\to coker(A)\to0.
 \eeq
The support of the module $coker(A)$ is the annihilator-of-cokernel ideal:
\beq
ann.coker(A)=ann\quotient{R^{\oplus m}}{Im(A)}=\{f\in R|\ f\cdot R^{\oplus m}\sseteq Im(A)\}
\sset R.
\eeq
 This ideal is $G_{lr}$-invariant and refines the ideal $I_m(A)$.

  The annihilator-of-cokernel is a rather delicate invariant but it is controlled by
 the ideals $I_j(A)$ as follows:
\beq\label{Eq.Background.ann.coker.in.terms.of.Fitting.ideals}\ber
\forall\ j<m:\ ann.coker(A)\cdot I_j(A)\sseteq I_{j+1}(A),\ {\rm and}\
ann.coker(A)^m\sseteq I_m(A)\sseteq ann.coker(A)\sseteq\sqrt{I_m(A)}\\
  \text{(proposition 20.7 in \cite{Eisenbud})}
\\
\text{If $m=n$ and $det(A)\in R$ is not a zero divisor, then $ann.coker(A)=I_m(A):I_{m-1}(A)$.}
\\
\text{If $m<n$ and $height(I_{m})=(n-m+1)$ then $ann.coker(A)=I_m(A)$, \cite[exercise 20.6]{Eisenbud}}
\eer\eeq
In particular, for one-row matrices, $m=1$, or when $I_m(A)$ is a radical ideal, $I_m(A)=ann.coker(A)$.

We use also the following properties:
\bprop\label{Thm.Background.Ann.Coker.Properties}
1. (Block-diagonal case) $ann.coker(A\oplus B)=ann.coker(A)\cap ann.coker(B)$.
\\2. If $A$ is a square matrix and $det(A)$ is not a zero divisor then $ann.coker(A)=ann.coker(A^T)$.
\\2'. If moreover $R$ is a unique factorization domain (UFD) then $ann.coker(A)$ is a principal ideal.
\eprop
Here statement {\bf 1} is immediate. Statement {\bf 2} follows directly from $ann.coker(A)=I_m(A):I_{m-1}(A)$ of equation (\ref{Eq.Background.ann.coker.in.terms.of.Fitting.ideals}).

For the statement {\bf 2'} note that for a square matrix the height of $ann.coker(A)\sset R$ is one. If $R$ is UFD then $ann.coker(A)$ is generated by just one element.

\subsection{The generalization of the annihilator of cokernel}\label{Sec.Background.Generaliz.ann.coker}
The ideal $ann.coker(A)$ is a `partially reduced' version of the ideal of maximal minors $I_m(A)$. (Alternatively, the annihilator of a module, $ann(M)$, is a refinement of the minimal Fitting ideal of that module, $Fitt_0(M)$.)
More generally, the counterparts of the ideals $\{I_j(A)\}$ (or the Fitting ideals $\{Fitt_{m-j}(M)\}$) are described in \cite{Buchsbaum-Eisenbud}, see also \cite[exercise 20.9]{Eisenbud}.
 We recall briefly the definition and the main properties.

Fix a morphism of free $R$-modules, $E\stackrel{\phi}{\to}F$, here $rank(F)=m$. For each $1\le j\le m$ define the associated morphism $E\otimes\wedgel^{j-1}F\stackrel{\phi_j}{\to}\wedgel^j F$ by $a\otimes w\to \phi(a)\wedge w$.
\bed
$ann.coker_j(\phi):=\Bigg\{\bM ann.coker(\phi_{m+1-j}),&\rm{for}\ 1\le j\le m\\
R,&j\le 0\\\{0\},& j>m\eM$.
\eed

\begin{Properties}\label{Thm.Background.Ann.Coker.Generalization}
1.  $ann.coker_j(\phi)=ann\Big(\wedgel^{m+1-j}coker(\phi)\Big)$, in particular the ideal is fully determined by the module $coker(\phi)=\quotients{F}{\phi(E)}$.
\\2. The ideals $ann.coker_j(\phi)$ refine the determinantal ideals, in the following sense:
\bee[i.]
\item
$ann.coker(\phi)=ann.coker_m(\phi)\sseteq\cdots\sseteq ann.coker_j(\phi)\sseteq ann.coker_{j-1}(\phi)\sseteq\cdots\sseteq
 ann.coker_1(\phi)=I_1(\phi)$.
\item For any $i>j\ge0$ holds: $ann.coker_i(\phi)\cdot I_j(\phi)\sseteq I_{j+1}(\phi)$.
\item $ann.coker_j(\phi)\supseteq I_j(\phi)\supseteq (ann.coker_j(\phi))^{j}$.
\item  $I_j(\phi)\sseteq ann.coker_j(\phi)\sseteq I_j(\phi):I_{j-1}(\phi)$. In particular, if $I_{j-1}(\phi)$ has a non-zero divisor modulo $I_j(\phi)$ then $ann.coker_j(\phi)=I_j(\phi):I_{j-1}(\phi)$.
\eee
3.i. Suppose the map splits block-diagonally, i.e. $E_1\oplus E_2\stackrel{\la\oplus\de}{\to}F_1\oplus F_2$. Suppose moreover $\la$ is invertible  (thus in particular $rank(E_1)=rank(F_1)$). Then
$ann.coker_j(\phi)=ann.coker_{j-rank(F_1)}(\de)$.
\\3.ii. If $A=diag(\la_1,\dots,\la_m)\in Mat_{m\times m}(R)$ and $(\la_1)\supseteq(\la_2)\supseteq\cdots\supseteq(\la_m)$ then $ann.coker_j(\phi)=(\la_k)$.
\\4. The ideals $ann.coker_j(\phi)$ are functorial under localizations, i.e. $ann.coker_j(\phi)_\cp=ann.coker_j(\phi_\cp)$ for any prime $\cp$.
\\5. Suppose $rank(Im(\phi))<r$, then $ann_j(\phi)=\{0\}$ for $j\ge r$.
\end{Properties}
Some remarks/explanations are needed here.
\\1. Fix some bases of $E,F$, so that $\phi$ is presented by a matrix $A$.  Then
 $ann.coker_j(\phi)$ is invariant under $G_{lr}$-action on $A$. Similarly, given  $A\in Mat_{m\times n_1}(R)$ and $B\in Mat_{m\times n_2}(R)$. If $Span(Columns(A))=Span(Columns(B))$ then $ann.coker_j(A)=ann.coker_j(B)$.     If $n_1=n_2$ and  $Span(Rows(A))=Span(Rows(B))$ then $ann.coker_j(A)=ann.coker_j(B)$.
\\2.i. This sequence of inclusions and the equalities are immediate.
\\2.ii. and 2.iii see  \cite[exercise 20.9]{Eisenbud} and \cite[exercise 20.10]{Eisenbud}. For 2.iv see corollary 1.4. of \cite{Buchsbaum-Eisenbud}.
\\3.i. In this case $coker(\phi)\approx coker(\de)$, now use part 1.
\\3.ii. Follows by explicit check.
\\4. Follows straight from $ann(M_\cp)=ann(M)_\cp$.
\\5. If $rank(Im(\phi))<r$ then $rank(Im(\phi_{m+1-r}))< rank(\wedgel^{m+1-r}F)$. But then $ann.coker(\phi_{m+1-r})=\{0\}$.

\subsection{The singular locus of an ideal}\label{Sec.Background.Sing(J)}
This is defined in \S\ref{Sec.Results.Sing(J)}, here we collect its basic properties.



Though the definition involved various choices of the bases, $Sing_r(J)$ depends on the ideal $J$ only, behaves well under localizations, sometimes equals $J$ and has other good properties.
\bel\label{Thm.Background.Sing(J).Properties}
Let $J\sset R$ be an ideal of expected height $r\le dim(R)$. (For the statements 3,4,5 we assume $R$ Noetherian.)
\\1. The ideal $Sing_r(J)$ does not depend on the choice of the generators of $J$ or $Der_\k(R)$.
\\2. $Sing_r(J)\supseteq J+ann_r\quotients{R^{\oplus k}}{Der_\k(R)(\{f_i\})}$, and the inclusion can be proper.
\\3. $\sqrt{Sing_r(J)}=\sqrt{J+I_r(Der_\k(R)(\{f_i\}))}$
\\4. If either $height(J)<r$ or $rank(Der_\k(R))<r$ then $Sing_r(J)=J$.
\\5. For any prime ideal $\cp\supset\cm$ holds: $Sing_r(J)_\cp=Sing_r(J_\cp)$. If $\cp\ssetneq\cm$ then
$Sing^{(\cm)}_r(J)_\cp=Sing_r(J_\cp)$.
\\6.i. For any matrix $A\in \Mat$ holds: $I_{r+1}(A)\sseteq Sing_{(m-r)(n-r)}(I_{r+1}(A))\sseteq I_r(A)$.
\\6.ii. For $A\in Mat^{sym}_{m\times m}(R)$ holds:
 $I_{r+1}(A)\sseteq Sing_{\bin{m+1-r}{2}}(I_{r+1}(A))\sseteq I_r(A)$.
\eel
\bpr
1.
Suppose $\uf$, $\ug$ are two tuples of generators of $J$. The two tuples are related by $\uf=U\ug$, $\ug= V\uf$, where $U,V$ are some matrices over $R$ (rectangular, not necessarily square).

Then the matrices in equation  \eqref{Eq.Definition.of.Sing(J)}, for $\uf$, $\ug$, are related by the left-right action by matrices (over $R$). Hence we get: $Sing^{(\uf)}(J)\sseteq Sing^{(\ug)}(J)\sseteq Sing^{(\uf)}(J)$, and thus $Sing^{(\uf)}(J)=Sing^{(\ug)}(J)$.

The independence of the choice of generators of $Der_\k(R)$ is even simpler, in this case one should apply only the right multiplications of the matrices.

2.  The inclusion is obvious and the following example shows the inequality. Take $J=(x^7+y^8,x^8+y^9)\sset \k[[x,y]]$. The height of this ideal is two, as expected. We claim that $Sing_2(J)\supseteq(x^7,y^8)$.
Indeed:
\beq
Sing_2(J)=ann.coker\bpm x^7+y^8& x^8+y^9&0&0&7x^6&8y^7
\\
0&0&x^7+y^8&x^8+y^9&8x^7&9y^8
\epm
\eeq
Now, a combination of the first, third and fifth columns of this matrix gives $\bpm 7y^8\\8y^8\epm$. Together with the sixth column this gives $\bpm y^7&0\\0&y^8\epm$, which ensures: $y^8\in Sing_2(J)$.  Now we get $\bpm x^6\\0\epm,\bpm0\\x^7\epm$, thus $x^7\in Sing_2(J)$.

But
 $J+ann_r\quotients{R^{\oplus k}}{Der_\k(R)(\{f_i\})}=(x^7+y^8,x^8+y^9,56x^6y^8-64x^7y^7)=(x^7+y^8,xy^8+y^9,\cm^6\cdot y^{8})$. (The last transition uses the Gr\"{o}bner basis.) From here one sees that e.g. $Sing_2(J)\ni y^8\not\in J+ann_r\quotients{R^{\oplus k}}{Der_\k(R)(\{f_i\})}$.

3. By lemma  \ref{Thm.Background.Localization.2} it is enough to prove:
\beq
\text{ for any prime  ideal $\cp\ssetneq\cm$ there holds: $Sing_r(J)_\cp=R_\cp$ iff $\Big(J+I_r(Der_\k(R)(\{f_i\}))\Big)_\cp=R_\cp$.}
\eeq
 If $\cp\not\supseteq J$ then both sides are $R_\cp$, as both sides  contain $J$, and $J_\cp=R_\cp$ by lemma \ref{Thm.Background.Localization.1}. If $\cp\supseteq J$ then
\beq
\text{$Sing_r(J)_\cp=R_\cp$ iff
 $\Big(ann_r\quotients{R^{\oplus k}}{Der_\k(R)(\{f_i\})}\Big)_\cp=R_\cp$ iff
$I_r(Der_\k(R)(\{f_i\}))_\cp=R_\cp$.}
\eeq

4.
Fix some generators $f_1,\dots,f_N$ of $J$ and present the singular locus in the form:
\beq
Sing_r(J)=ann_r\quotient{R^{\oplus N}}{J\cdot R^{\oplus N}+Der_\k(R)\bpm f_1\\\dots\\f_N\epm}=ann_r\quotient{\Big(\quotients{R}{J}\Big)^{\oplus N}}{\quotients{R}{J}\otimes  Der_\k(R)\bpm f_1\\\dots\\f_N\epm}.
\eeq

{\bf i.} Suppose $height(J)<r$, then for any $r$-tuple of elements $f_{i_1},\dots,f_{i_r}\in J$ there exists a non-Koszul relation, $\suml^r_{j=1}a_{i_j}f_{i_j}=0$, i.e. not all the coefficients $\{a_{i_j}\}$ belong to $J$. Thus such a relation gives a non-trivial relation over $\quotients{R}{J}$ as well. And the later leads to a non-trivial relation:
  $\suml^r_{j=1}a_{i_j}\cD f_{i_j}=0\in \quotients{R}{J}$, for any derivation $\cD\in Der_\k(R)$.  Therefore any $r$ rows of the matrix $\quotients{R}{J}\otimes  Der_\k(R)(\{f_i\})$ are linearly dependent over $\quotients{R}{J}$. In other words, the rank of module
 $\quotients{R}{J}\otimes  Der_\k(R)(\{f_i\})$ is less than $r$. And thus, by the property
 6 of \ref{Thm.Background.Ann.Coker.Generalization} we have:
\beq
ann_r\quotient{(R/J)^{\oplus k}}{(R/J)\otimes Der_\k(R)(\{f_i\})}=ann_r\Big((\quotients{R}{J})^{\oplus k}\Big)=J.
 \eeq

{\bf ii.}  Suppose the rank of $Der_\k(R)$ is less than $r$. We use the base change:
$\quotients{R}{J}\otimes ann^R_r(\dots)=ann^{\quotients{R}{J}}_r(\dots)$. And now use property 6 of  \ref{Thm.Background.Ann.Coker.Generalization} to get $ann_r(\quotients{R}{J}\otimes\dots)=\{0\}$. Finally, we have:
$\quotients{R}{J}\otimes Sing_r(J)=\{0\}$ and of course $Sing_r(J)\supseteq J$, thus $Sing_r(J)= J$.

5.
The equality $Sing_r(J)_\cp=Sing_r(J_\cp)$ holds because the annihilator is functorial, $ann_r(M_\cp)=ann_r(M)_\cp$, and the module of derivations as well: $Der_\k(R)_\cp=Der_\k(R_\cp)$.

If $\cp\ssetneq\cm$ then $Der_\k(R,\cm)_\cp\supseteq (\cm\cdot Der_k(R))_\cp=Der_k(R)_\cp=Der_k(R_\cp)$.
Therefore $Sing^{(\cm)}_r(J)_\cp=Sing_r(J_\cp)$.


6.
 The inclusion $I_{r+1}(A)\sseteq Sing_{(m-r)(n-r)}(I_{r+1}(A))$ holds by part 2. For the inclusion
$Sing(I_{r+1}(A))\sseteq I_{r}(A)$ note that $Der_\k(R)(I_{r+1}(A))\sseteq I_r(A)$. (Expand the $(r+1)\times(r+1)$ determinants in terms of $r\times r$ determinants.) Therefore
\beq
ann_{(m-r)(n-r)}\quotient{R^{\oplus N}}{I_{r+1}(A)\cdot R^{\oplus N}+Der_\k(R)\bpm f_1\\\dots\\ f_N\epm}\sseteq
ann_{(m-r)(n-r)}\quotient{R^{\oplus N}}{I_{r+1}(A)\cdot R^{\oplus N}+I_{r}(A)\cdot R^{\oplus N}}=I_r(A).
\eeq
The symmetric case is similar.
\epr

\subsection{Tangent spaces to the orbits}\label{Sec.Background.Tangent.Spaces}
The tangent spaces for various group actions are written in \cite[\S3.5]{Belitski-Kerner}.

\bex\label{Ex.Tangent.Spaces.to.the.Actions}
For the group-actions of example \ref{Ex.Intro.Typical.Groups} the tangent spaces are $R$-modules, we recall their presentation.
\li $G_{lr}:\ A\to UAV$. Here $T_{(G_{lr}A,A)}=\Span_R\{uA,Av\}_{(u,v)\in Mat_{m\times m}(R)\times Mat_{n\times n}(R)}$. Similarly
 for $G_l$ and $G_r$.
\li $Aut_\k(R)$. Any automorphism of the ring, $\phi\circlearrowright R$, acts on the matrices entry-wise: $\{a_{ij}\}\to \{\phi(a_{ij})\}$. The tangent space to the orbit is obtained by applying the tangent space $T_{(Aut_\k(R),Id)}$ to $A$. We have $T_{(Aut_\k(R),Id)}=Der_\k(R,\cm)$,  the module of those derivations of $R$ that send $\cm$ into itself.
 (Here we have only the submodule $Der_\k(R,\cm)\sseteq Der_\k(R)$ because the automorphisms of a local ring preserve the origin of $Spec(R)$.) Therefore
\beq
 T_{(Aut_\k(R)(A),A)}=Der_\k(R,\cm)(A)=\Span_R\{\cD(A)\}_{\cD\in Der_\k(R,\cm)}.
\eeq
For a regular subring of $\k[[\ux]]$ holds: $Der_\k(R,\cm)=\cm\cdot Der_\k(R)$, where the module  $Der_\k(R)$
 is generated by the first order partial derivatives $\{\di_j\}$.
\li $\cG_{lr}:\ A\to U\phi(A)V$. Here
$T_{(\cG_{lr}A,A)}=\Span_R\{uA,Av,\cD(A)\}_{(u,v,\cD)\in Mat_{m\times m}(R)\times Mat_{n\times n}(R)\times Der_\k(R,\cm)}$.
Similarly  for $\cG_l$ and $\cG_r$.
\li $\cG_{congr}:\ A\to U\phi(A)U^t$. Here
$T_{(\cG_{congr}A,A)}=\Span_R\{uA+Au^t,\cD(A)\}_{(u,\cD)\in Mat_{m\times m}(R)\times Der_\k(R,\cm)}$.
\eex

\bex
To compute/bound the order of determinacy we need also the tangent spaces to the orbits of $G^{(J)}$ for some ideal $J\sset R$. These are obtained from the previous example by using: $T_{(G^{(J)}_{lr}A,A)}=J\cdot T_{(G_{lr}A,A)}$, and replacing $Der_\k(R,\cm)$ by $Der_\k(R,\cm\cdot J)$.
\eex

\subsection{Invariance of the annihilator}\label{Sec.Background.Invariance.Annihilator}
The element $h=(U,V,\phi)\in\cG_{lr}$ acts on $R$ by $f\to\phi^*(f)$. We use the sloppy notation $h^*(f)$ and $h^*(J)$ for an ideal $J\sset R$.
Suppose $h\in \G_{lr}$ acts on $\Si$, i.e. it sends the germ $(\Si,A)$ to the germ $(\Si,hA)$.
\bprop\label{Thm.Annihilator.Invariance.Group.Action}
 Let $R$ be a local ring. Suppose the tangent spaces $T_{(GA,A)}$, $T_{(\Si,A)}$ are $R$-modules.
  Suppose $h\in\cG_{lr}$ acts on $\Si$ and also commutes with the $G$-action on $A$, i.e. $hGA=GhA$ .
 Then $h^*(ann(T^1_{(\Si,G,A)}))=ann(T^1_{(\Si,G,h(A))})$.
\eprop
 In particular $ord^\Si_G(A)=ord^\Si_G(hA)$.
\bpr
Consider $h$ as a $\k$-linear automorphism of $\Mat$. It induces the isomorphism of the tangent spaces:
\\$T_{(\Mat,A)}\stackrel{h_*}{\isom}T_{(\Mat,hA)}$. As $h$ acts on $\Si$, the tangent isomorphism restricts to
 $T_{(\Si,A)}\stackrel{h_*}{\isom}T_{(\Si,hA)}$. Further, the restriction $(GA,A)\stackrel{h}{\isom}(hGA,hA)=(GhA,hA)$ induces the isomorphism $T_{(GA,A)}\stackrel{h_*}{\isom}T_{(GhA,hA)}$.

If $h\in G_{lr}$ then the map $h_*$ is $R$-linear. If $h\in\cG_{lr}$ then the map is compatible with  $R$-multiplication:
\beq
h_*(f\cdot T_{(\Si,A)})=h^*(f)\cdot h_*(T_{(\Si,A)}),\quad  h_*(f\cdot T_{(GA,A)})=h^*(f)\cdot h_*(T_{(GA,A)}).
\eeq
 Altogether:
$h_*\big(f\cdot T^1_{(\Si,G,A)}\big)=h^*(f)h_*\big(T^1_{(\Si,G,A)}\big)$.
 Thus, if $f\in ann\big(T^1_{(\Si,G,A)}\big)$ then
 $h^*(f)\in ann\big(T^1_{(\Si,G,hA)}\big)$, i.e.
 $h^*\big(ann T^1_{(\Si,G,A)}\big)\sseteq ann(T^1_{(\Si,G,hA)})$.
As $h$ is invertible, we get the inverse inclusion as well.\epr

\bex
$\bullet$ In the simplest case, $h\in G$, we get the obvious property: the $G$-determinacy is constant along the $G$-orbit.
\li In many cases no choices of $h\in G$ helps,  e.g. $A$ has no nice canonical form under the $G$-action.
Then one extends $G$ by its normalizer, as in the proposition. For example, we use the following normal extensions: $G_l,G_r\triangleleft G_{lr}$, $\cG_r\triangleleft GL(m,\k)\times\cG_r$, $\cG_l\triangleleft \cG_l\times GL(n,\k)$. \li Note that $h\in\cG_{lr}\smin\{GL(m,\k)\times\cG_r\}$ does not in general normalize the $\cG_r$-action.
 Similarly for $h\in\cG_{lr}\smin\{\cG_l\times GL(n,\k)\}$.
\eex

\subsection{The relevant approximation properties}\label{Sec.Background.Approximation.Properties}
 (This continues \S\ref{Sec.Results.Relevant.Approx}.) Fix some group action  $\cG_{lr}\supseteq G\circlearrowright\Mat$.
The completion of the ring, $R\to\hR$, induces the map $\Mat\to Mat_{m\times n}(\hR)$ and accordingly the completion map of the groups: $\cG_{lr}\supseteq G\to \hG\sseteq \widehat{\cG_{lr}}$, see \cite[\S3.2]{Belitski-Kerner} for the details.
  We often need the following approximation property:
\beq\label{Eq.Approximation.Property.for.Groups}\ber
\text{\em Given the equation $gA=A+B$ for $g\in G$, suppose there exists a formal solution, $\hat{g}\hA=\hA+\hB$, $\hat{g}\in\hG$.}
\\ \text{\em  Then there exists an ordinary solution, $gA=A+B$, $g\in G$.}
\eer\eeq
This property restricts the possible rings, the particular condition on $R$ depends on the type of the equations, see \cite[\S5]{Belitski-Kerner} for the details/proofs.
\bee[i.]
\item If $G\sseteq G_{lr}$ is defined by $R$-linear equations  and the condition $gA=A+B$ can be written as a system of linear equations on $g=(U,V)$ then the property (\ref{Eq.Approximation.Property.for.Groups}) holds for $G$ over arbitrary Noetherian local ring $R$.

In the non-Noetherian case the property holds if
the completion is surjective, $R\twoheadrightarrow\hR$ and $ann.coker(A)\supseteq\cm^\infty$.
 Note that by Borel's lemma: $R=C^\infty(\R^p,0)\twoheadrightarrow\R[\![\ux]\!]$, \cite[pg. 284, exercise 12]{Rudin-book}. Further, the condition $ann.coker(A)\supseteq\cm^\infty$ can checked
 using the property
\begin{multline}\label{Eq.Lojasiewicz.Ineq}
\text{If $f\in C^\infty(\R^p,0)$ satisfies $|f|\ge C|\ux|^\de$, for some $C>0$, $\de>0$,}
 \\\text{then $f$ divides any function flat at the origin, and thus $(f)\supseteq\cm^\infty$.}
\end{multline}
\item If $G\sseteq \cG_{lr}$ is defined by polynomial/analytic equations and the condition $gA=A+B$ can be written as a system of polynomial/analytic
equations on $g=(U,V,\phi)$,
then the property (\ref{Eq.Approximation.Property.for.Groups}) holds over any Henselian Noetherian ring. When the group involves
 the coordinate changes and $\phi\neq Id$ then the assumption implies that $A$ is a matrix of polynomials/analytic functions.
\item If the defining equations of $G\sseteq\cG_{lr}$ are not polynomial/analytic
 or the condition $gA=A+B$ cannot be written polynomially/analytically in the entries of $g=(U,V,\phi)$, then the approximation  (\ref{Eq.Approximation.Property.for.Groups}) holds at least when the ring $R$ is a quotient of a Weierstrass system, \cite[Example 2.13]{Hauser-Rond}. The simplest examples of  Weierstrass system are:
 the formal power series, $\k[\![\ux]\!]$; the algebraic power series, $\k<\ux>$; the analytic power series, $\k\{\ux\}$; Gevrey power series. (In the last two cases $\k$ is a normed field.)
\eee

\subsection{The completion of $Aut_\k(R)$ and the condition $\widehat{T_{(Aut_\k(R),Id)}}=T_{(\widehat{Aut_\k(R)},Id)}$}\label{Sec.Background.Aut(R)}
Let $R$ be a local ring.
The filtration $\{\cm^j\}$ on $R$ induces the filtration of $Aut_\k(R)$ by the subgroups
$Aut^{(\cm^j)}_\k(R):=\{\phi\in Aut_\k(R)|\ \phi\equiv Id\ mod\ \cm^j\}.$

These subgroups are normal, $Aut^{(\cm^j)}_\k(R)\vartriangleleft Aut_\k(R)$, and one often considers the completion \wrt this filtration, $\widehat{Aut_\k(R)}:=\liml_{\leftarrow}\quotients{Aut_\k(R)}{Aut^{(\cm^j)}_\k(R)}$.

This group completion acts naturally on the ring completion, $\hR:=\liml_{\leftarrow}\quotients{R}{\cm^j}$, the action is defined as follows.
For any
$\hf\in \hR$ choose a representing sequence, $\{f_j\in \quotients{R}{\cm^j}\}$.
For any $\hphi\in \widehat{Aut_\k(R)}$ choose a representing sequence, $\{\phi_j\in \quotients{Aut_\k(R)}{Aut^{(\cm^j)}_\k(R)}\}$.
Define
\beq
\hphi(\hf)=\{\phi_j\}(\{f_i\}):=\{\phi_j(f_j)\}.
\eeq
(By the direct check: the action $\hphi(\hf)$ does not depend on the choice of representatives for $\hphi$, $\hf$.)

This action induces the group homomorphism $\widehat{Aut_\k(R)}\to Aut_\k(\hR)$. This map is injective but in general not surjective.
\bex
Two reasons for the possible non-surjectivity are: flat functions and/or non-Henselianity of $R$.
\bee
\item Fix a smooth function-germ in one variable, $\tau\in C^\infty(\R^1,0)$, satisfying: $\tau$ is flat at $0\in \R$ and $\tau|_{\R^1\smin\{0\}}>0$. Take the ring $R=\quotients{C^\infty(\R^p,0)}{(\tau)}$, then $\hR=\R[[\ux]]$.
 We have: $Aut_\k(R)=\{\phi\circlearrowright C^\infty(\R^p,0)|\ \phi(\tau)=\tau\}$. Therefore $\widehat{Aut_\k(R)}\ssetneq Aut_\k(\hR)$.
\item Consider the localization of affine ring, $R=\Big(\quotients{\k[x,y]}{(y^2-x^2-x^3)}\Big)_\cm$. Then
$\hR\approx\quotients{\k[[\tx,y]]}{(y^2-\tx^2)}$, here $\tx=x\sqrt{1+x}$. Therefore $Aut_\k(\hR)$ contains the permutation of axes, $\tx\leftrightarrow y$, while $\widehat{Aut_\k(R)}$ does not contain such an element.
\eee
\eex

\bprop\label{Thm.Background.T(Aut(R)).and.completion}
Suppose
\bee[i.]
\item either $R=\quotients{S}{I}$, where $S\sseteq\k[[\ux]]$ is a regular Henselian subring, closed under the action of partial derivatives, $\di_i(S)\sseteq S$, and whose completion is $\k[[\ux]]$.
(For example $R$ is one of $\quotients{\k[[\ux]]}{I}$, $\quotients{\k\{\ux\}}{I}$, $\quotients{\k<\ux>}{I}$, note that the derivative of an algebraic series is algebraic.)
\item or $R=\quotients{C^\infty(\R^p,0)}{I}$, where the ideal $I$ is finitely generated by algebraic power series.
\eee
1. Then $\widehat{Aut_\k(R)}=Aut_\k(\hR)$ and $\widehat{Der_\k(R)}=Der_\k(\hR)$.
\\2. In particular, for these rings holds: $\widehat{T_{(Aut_\k(R),\one)}}=T_{(\widehat{Aut_\k(R)},\one)}$.
\eprop
(Here $\k<\ux>$ is the ring of algebraic power series.)
\bpr
1. $\bullet$ {\em The case of $R=\quotients{S}{I}$.}
By the assumption $\hR=\quotients{\k[[\ux]]}{\hI}$ and every element of $Aut_\k(\hR)$ comes from some $\hphi\in Aut_\k(\k[[\ux]])$ that satisfies $\hphi(\hI)=\hI$. Choose some  generators $\{g_i(\ux)\}$
 of  $I\sset S$, (a finite set!), then $\hphi$ is a tuple of power series that satisfy the system of implicit function equations:
\beq
\forall\ j:\quad g_j(\hphi(\ux))=\hphi(g_j(\ux))=\sum z_i g_i(\ux),\ z_i\in \hat{S}.
\eeq
But $S$ is Henselian, thus for any $N$ there exists a solution, $\phi_N$ in $S$, that satisfies these equations and coincides with $\hphi$ modulo $\cm^N$. Thus there exists a sequence of elements of $Aut_\k(S)$ that converges ($\cm$-adically) to $\hphi\in Aut_\k(\hS)$. As this sequence preserves $I\sset S$ we get an element in $\widehat{Aut_\k(R)}$ that goes to the element in $Aut_\k(\hR)$.

For the derivations, we first note that the derivations of $\k[[\ux]]$ are determined by their action on the coordinates $\ux$, i.e. $\cD(f)=\sum \cD(x_i)\di_i(f)$, here $\di$ are the ordinary partial derivatives. The same holds for the subring $S\sseteq\k[[\ux]]$.
Every element of $Der_\k(\hR)$ comes from a derivation $\widehat{\cD}\in Der_\k(\hS)$ that satisfies $\widehat{\cD}(\hI)\sseteq \hI$. Expand $\widehat{\cD}=\sum \hf_i\di_i$.
We get the system of equations, $\{\sum \hf_i\di_i(g_j)=\sum a_{jk}g_k\}$, with a solution $\hf_i,a_{jk}\in \k[[\ux]]$. Then, as before, this solution can be approximated by a solution in $S$. Which means:
 $\widehat{\cD}$ is the limit of some sequence of derivations, $\cD_N\in Der_\k(S)$, that satisfy this system of equations, i.e. each of them descends to $Der_\k(R)$. This gives a Cauchy sequence of elements $Der_\k(R)$ that converges to the element of $Der_\k(\hR)$.

$\bullet$ {\em The case of $R=\quotients{C^\infty(\R^p,0)}{I}$.}
Fix an automorphism of $\quotients{C^\infty(\R^p,0)}{I}$, it is presented by some $\phi\circlearrowright C^\infty(\R^p,0)$ that satisfies $\phi(I)=I$. In particular the action of $\phi$ on $\ux$ is fixed. Define the associated change-of-variables automorphism $\phi_{geom}\circlearrowright C^\infty(\R^p,0)$ by $\phi_{geom}(f(x))=f(\phi(x))$. We claim that $\phi_{geom}(I)=I$ and therefore $\phi_{geom}$ descends to an automorphism of $\quotients{C^\infty(\R^p,0)}{I}$.

Indeed, by its definition the automorphism $\phi\circ\phi^{-1}_{geom}$ acts as the identity on the polynomials. In particular, $\phi\circ\phi^{-1}_{geom}\in Aut^{\cm^\infty}_\R(C^\infty(\R^p,0))$.
 We claim that $\phi\circ\phi^{-1}_{geom}$ acts as the identity on the algebraic power series. Indeed, any algebraic power series $f$ satisfies a polynomial equation over $\k[\ux]$, therefore $\phi\circ\phi^{-1}_{geom}(f)$ satisfies the same equation. But $f-\phi\circ\phi^{-1}_{geom}(f)\in\cm^\infty$, therefore $f=\phi\circ\phi^{-1}_{geom}(f)$.
Therefore, as $I$ is generated by algebraic power series, $(\phi\circ\phi^{-1}_{geom})(I)=I$.

Therefore  $\phi_{geom}(I)=I$, hence $\phi_{geom}$ defines an automorphism of $\quotients{C^\infty(\R^p,0)}{I}$.
 Therefore
$\widehat{Aut_\R(R)}$ coincides with the completion of the subgroup of geometric automorphisms. But  the later completion is just $Aut_\k(\hR)$.

To any derivation $\cD\in Der_\k(R)$ associate the ``geometric derivation", $\cD_{geom}(f)=\sum \cD(x_i)\di_i(f)$. Then $\cD-\cD_{geom}\in Der_\k^{\cm^\infty}(R)$ and therefore $\widehat{Der_\k(R)}$ coincides with the completion of the subspace of geometric derivations. But the later completion is just $Der_\k(\hR)$.

2. This follows immediately from $T_{(Aut_\k(R),\one)}=Der_\k(R)$.
\epr

\subsection{The general linearization statement: transition from the annihilator of $T^1_{(\Si,G,A)}$ to finite  determinacy}\label{Sec.Background.BK}
Let $M$ be a finitely generated $R$-module with  a group action  $G\circlearrowright M$. Suppose the action preserves a subset $\Si\sseteq M$.
In \cite{Belitski-Kerner} we have reduced the study of determinacy to the understanding of the support/annihilator of the module $T^1_{(\Si,G,z)}$.  Below we formulate the results used in the current paper.

\

Fix an ideal $J\sset R$, this gives a filtration on the module, $\{J^q\cdot M\}$, and the corresponding filtration on the group, $\{G^{(J^q)}\}$, see example \ref{Ex.Intro.Typical.Groups}.

For this section we assume that the group actions $G\circlearrowright M$ and $G^{(J)}\circlearrowright M$ satisfy the following conditions:
\beq
\text{the germ $(G,\one)$ is formally smooth and the natural map $\widehat{T_{(G,\one)}}\to T_{(\widehat{G},\one)}$ is an isomorphism.}
\eeq
These are satisfied for various subgroups of $G_{lr}$, \cite[\S 3.8]{Belitski-Kerner}. If the ring automorphisms (coordinate changes) are involved then the condition $\widehat{T_{(G,\one)}}\isom T_{(\widehat{G},\one)}$ implies some restrictions on $R$, see \S\ref{Sec.Background.Aut(R)}, yet proposition \ref{Thm.Background.T(Aut(R)).and.completion} ensures this condition for our rings.

\bprop \cite[Corollary 2.5]{Belitski-Kerner}\label{Thm.Background.Reduction.to.annihilator}
 Suppose $T_{(\Si,z)}\sseteq T_{(M,z)}$ is a finitely generated submodule and for a (finitely-generated) ideal $J\sseteq \cm$ the filtration $\{(\sqrt{J})^q\cdot T_{(\Si,z)}\}_q$ is $G$-invariant.
\\1. Suppose $R$ has the relevant approximation property (\S\ref{Sec.Background.Approximation.Properties}) and $J\sseteq ann(T^1_{(\Si,G,z)})$. Then $\{z\}+ J\cdot\sqrt{J}\cdot T_{(\Si,z)}\sseteq Gz$.
\\1'.  If in addition $J\cdot T_{(\Si,z)}\sseteq T_{(G^{(\sqrt{J})}z,z)}$ then $\{z\}+J\cdot T_{(\Si,z)}\sseteq Gz$.
\\2. If $\{z\}+ J\cdot T_{(\Si,z)}\sseteq G^{(\sqrt{J})}z$ then $J\sseteq ann(T^1_{(\Si,G^{(\sqrt{J})},z)})$.
\eprop

This proposition reduces the determinacy question to the study of the module $T^1_{(\Si,G,A)}$ and its annihilator.

In particular, it allows to bound the order of determinacy as follow.
(The Loewy length of an ideal, $ll_R(..)$, is defined in \S\ref{Sec.Results.Ideals}.)
\bprop\label{Thm.Background.Ord.of.det.via.Loewy.length}\cite[Corollary 2.7]{Belitski-Kerner}
Suppose $\Si\sseteq\Mat$ is a free direct summand, i.e. $\Si\oplus\Si^\bot=\Mat$ for a free submodule $\Si^\bot\sset\Mat$.

1. Suppose $ann(T^1_{(\Si,G,A)})\supseteq\cm^\infty$ and $R$ has the relevant approximation property (\S\ref{Sec.Background.Approximation.Properties}). Then
\[
ll_R\Big(ann(T^1_{(\Si,G,A)})\Big)-1\le ord^\Si_G(A)\le
ll_R\Big(ann(T^1_{(\Si,G^{(\cm)},A)})\Big)-1.
\]

2. If $ann(T^1_{(\Si,G,A)})\not\supseteq\cm^\infty$, then $A$ is not infinitely-$(\Si,G)$-determined.
\eprop

In particular, $A\in \Mat$ is finitely determined iff $\cm^N T^1_{(\Si,G,A)}=\{0\}$
 for some $N<\infty$, alternatively: $\cm^N T_{(\Si,A)}\sseteq T_{(GA,A)}$.

\section{Proofs, corollaries and examples}\label{Sec.Proofs}
\subsection{The automorphisms of the ring, $Aut_\k(R)$}\label{Sec.Proofs.Aut(R)}
 \

{\em proof of Theorem \ref{Thm.Results.T1.for.Aut(R)}.}

Part 1 follows directly from the presentation of $T_{(Aut_\k(R)A,A)}$, \S\ref{Sec.Background.Tangent.Spaces}.

Part 2 follows now by proposition \ref{Thm.Background.Reduction.to.annihilator}.  To use this proposition  we note that the filtration $\{Mat_{m\times n}((\sqrt{J})^q)\}_q$ is invariant under the action of $Aut^{(\sqrt{J})}_\k(R)$. 

Part 3 now follows by proposition \ref{Thm.Background.Ord.of.det.via.Loewy.length}, note that $T_{(Aut^{(\cm)}_k(R),\one)}=Der_\k(R,\cm^2)$.
\epr

\bex\label{Ex.Proofs.Relative.Determinacy.Aut(R)}  (Continuing example \ref{Ex.Results.Fin.Determ.Aut(R)}.)
Suppose $I_1(Jac^{(\cm)}(f))$ contains no finite power of $\cm$, thus $f$ is not finitely-$Aut_\k(R)$-determined. Take the saturation,
\beq
Sat_\cm(I_1(Jac^{(\cm)}(f)))=\suml_j I_1(Jac^{(\cm)}(f)):\cm^j=\suml_j I_1(Jac(f)):\cm^j.
\eeq
 Suppose $J\sseteq Sat_\cm(I_1(Jac^{(\cm)}(f)))$.
 Then $\Si=J\sqrt{J}\sset R$ can serve as a module of admissible deformations in the following sense.
 Consider those ring automorphisms/coordinate changes that preserve $\sqrt{J}$, namely  $Aut^{(0)}_\k(R)=\{\phi\in Aut_\k(R)|\  \phi(\sqrt{J})=\sqrt{J}\}$. The tangent space to the orbit under this subgroup is
 \beq
 T_{(Aut^{(0)}_\k(R),Id)}=\{D\in Der_\k(R,\cm)|\ D(\sqrt{J})\sseteq \sqrt{J}\}\supseteq \sqrt{J}\cdot T_{(Aut_\k(R),Id)}.
 \eeq
Therefore we have:
\beq
ann\Big( T^1_{(\Si,Aut^{(0)}_\k(R),f)}\Big)=T_{(Aut^{(0)}_\k(R),Id)}(f): (J\sqrt{J})\supseteq \Big(\sqrt{J}\cdot T_{(Aut_\k(R),Id)}(f)\Big): (J\sqrt{J})\supseteq T_{(Aut_\k(R),Id)}(f): (J).
\eeq
 Here $T_{(Aut^{(0)}_\k(R),Id)}(f)$, $T_{(Aut_\k(R),Id)}(f)$ are considered as ideals in $R$.
 Note that now $ann\Big( T^1_{(\Si,Aut^{(0)}_\k(R),f)}\Big)$ contains some finite power of the maximal ideal, $\cm^N$. Thus the order of determinacy for these deformations satisfies:
\beq
ll_R\Big(T_{(Aut^{(0)}_\k(R),Id)}(f): (J\sqrt{J})\Big)-1\le ord^\Si_{Aut^{(0)}_\k(R)}(f)\le ll_R\Big(\cm\cdot T_{(Aut^{(0)}_\k(R),Id)}(f): (J\sqrt{J})\Big)-1.
\eeq
In the simplest case, suppose $R=\k[[x,y,z]]$, with $\k$ algebraically closed.
 Suppose the singular locus of $f$ is a smooth curve-germ (as a set), then after a change of coordinates it is   defined by the ideal $J=(x,y)\sset R$.
 Suppose the generic multiplicity of $f$ along the singular locus is $p=2$. We get:
  $J=Sat_\cm( I_1(Jac(f))$, so the space of admissible deformations is $J^2$.
 Then $ll_R(I)-1\le ord^\Si_{Aut_J(R)}(f)\le ll_R(I)$, where $I=\cm(J\di_1 f,\dots,J\di_k f,\di_{k+1}f,\dots,\di_p f):J^2$.
Compare to \cite{Siersma83}, \cite{Pellikaan}, \cite{Jiang}, see also \cite[Theorem 4]{Sun-Wilson}, \cite[Theorem 3.5]{Grandjean00}, \cite[Theorem 2.1]{Thilliez}.
\eex

{\em Proof of proposition \ref{Thm.Finite.Determinacy.Coord.Change}}.

$\Lleftarrow$ is obvious.

$\Rrightarrow$
By theorem \ref{Thm.Results.T1.for.Aut(R)} it is enough to check whether the ideal $ann.coker(Jac^{(\cm)}(A))$ contains some power of the maximal ideal. Here instead of $Jac^{(\cm)}(A)$ we can consider the matrix $Jac(A)$.

By the assumption, $Der_\k(R)=R<\di_1,\dots,\di_p>$ is a free $R$-module of rank $p$. Thus $\{\cD_\al a_{ij}\}_{\cD_\al\in Der_\k(R)}$ is a $n\times p$ matrix with the rows $(\di_1 a_{ij},\dots,\di_p a_{ij})$.

If $n>dim(R)$ then $ann.coker(\{\cD_\al a_{ij}\}_{\cD_\al\in Der_\k(R)})=0$, as there are more rows than columns.

For $n\le dim(R)$ the ideal $ann.coker(\{\cD_\al a_{ij}\}_{\cD_\al\in Der_\k(R)})$ is either non-proper, i.e. $R$, or of height at most $(dim(R)+1-n)$, which is less than $dim(R)$. Thus this ideal contains a power of the maximal ideal iff $ann.coker(\dots)=R$. Which means: at least one of the maximal minors is invertible. But the later means precisely that
 the entries of $A$ form a sequence of generators of $\cm$ (over $R$). Then the stability follows.
\epr
We remark that if $R$ is not regular then $A^{-1}(0)$ is not necessarily reduced or has an isolated singularity, see
example \ref{Example.Non.CM.ring} and remark \ref{Ex.Fin.Det.Nonisolated.Sing.for.non-regular.ring}.

\subsection{The $\cG_{lr}$-action}\label{Sec.Proofs.cG_lr}

\subsubsection{The set-theoretic support of $T^1_{(\Si,\cG_{lr},A)}$ for the geometric rings}\label{Sec.Proof.Supp(T1).Glr.set.theoretic}
For the sake of exposition we first prove the ``set-theoretic version" of part 4 of theorem \ref{Thm.Results.annT1.Grl.coord.changes} in the particular case:  for $R=\C\{\ux\}$.

Let $A\in\Mat$ and fix a small neighborhood  of the origin, $\cU\sset\C^p$, in which all the entries of $A$ converge. This neighborhood is stratified by the determinantal loci: $\cU=\coprod\limits^{m}_{r=0} \Big(V(I_{r+1}(A))\smin V(I_j(A))\Big)$. The formula for the support uses a refined stratification.
\bprop\label{Thm.Results.Support.Set.Theoretic.Glr.coord.change}
Let $R=\C\{\ux\}$.
The module $T^1_{(\Si,\cG_{lr},A)}$ is supported (set-theoretically) on the union of the complements,
\[\Supp(T^1_{(\Si,\cG_{lr},A)})=\coprod\limits^{m-1}_{r=0}\overline{V\Big(Sing_{(m-r)(n-r)}(I_{r+1}(A))\Big)\smin V\Big(I_r(A)\Big)}.
\]
\eprop
\bpr {\bf Step 1.}
As the statement is set-theoretic it is enough to check for any point $pt\in \cU$ that $pt$ belongs to the set on the right hand side iff it belongs to the set on the left hand side.

Fix $pt\in \cU$ and let $j$ satisfy $pt\in V(I_{j+1}(A))\smin V(I_j(A))$. Such $j$  exists and satisfies $m\ge j\ge0$ because $I_0(A)=R$, $V(I_0(A))=\varnothing$, and $I_{m+1}(A)=\{0\}$, $V(I_{m+1}(A))=\cU$.

We claim that $pt$ belongs to the union on the right hand side of the statement iff it belongs  to the particular term: $V\Big(Sing_{(m-j)(n-j)}(I_{j+1}(A))\Big)\smin V\Big(I_j(A)\Big)$. Indeed, $pt$ cannot belong to any term with $r<j$, as the erased part, $V\Big(I_r(A)\Big)$, contains $V\Big(I_{j+1}(A)\Big)$. And $pt$ cannot belong to any term with $r>j$ as $pt\not\in V\Big(I_j(A)\Big)$.

Thus we are to prove:
\begin{multline}\label{Eq.inside.proof}
\text{a point $pt\in V\Big(I_{j+1}(A)\Big)\smin V\Big(I_j(A)\Big)$ belongs to $\Supp(T^1_{(\Si,\cG_{lr},A)})$}
\\\text{iff it belongs to $V\Big(Sing_{(m-j)(n-j)}(I_{j+1}(A))\Big)\smin V\Big(I_j(A)\Big)$.}
\end{multline}
 Or:
\beq
\forall\ j:\quad \Supp(T^1_{(\Si,\cG_{lr},A)})\cap \Big(V\big(I_{j+1}(A)\big)\smin V\big(I_j(A)\big)\Big)=
V\Big(Sing_{(m-j)(n-j)}(I_{+1}(A))\Big)\smin V\Big(I_j(A)\Big).
\eeq
(If for some $j$ there holds $V\big(I_{j+1}(A)\big)=V\big(I_j(A)\big)$ then the condition is satisfied trivially.)

{\bf Step 2.} If $pt\in V\big(I_{j+1}(A)\big)\smin V\big(I_j(A)\big)$ then after a $G_{lr}$-transformation we can assume that  locally near $pt\in\cU$ the matrix is block-diagonal: $A=\one_{j\times j}\oplus\tA$. Here the matrix $\tA\in Mat_{(m-j)\times(n-j)}(R_{pt})$ vanishes at $pt$ and $R_{pt}$ denotes the local ring of $(\cU,pt)$.
(By proposition \ref{Thm.Annihilator.Invariance.Group.Action} the support is preserved under the $G_{lr}$-action.)

Therefore the tangent space decomposes into the direct sum, $T_{(\cG_{lr}A,A)}\approx Mat_{j\times n}(R_{pt})\oplus Mat_{(m-j)\times j}(R_{pt})\oplus T_{(\widetilde{\cG_{lr}}\tA,\tA)}$. Accordingly, near $pt$ there holds: $T^1_{(\Si,\cG_{lr},A)}\approx T^1_{(\tilde\Si,\widetilde\cG_{lr},\tA)}$, where $\tilde\Si=Mat_{(m-j)\times(n-j)}(R_{pt})$ and $\widetilde\cG_{lr}$ is the corresponding group.

Finally, $pt\in \Supp(T^1_{(\tilde\Si,\widetilde\cG_{lr},\tA)})$ iff $T^1_{(\tilde\Si,\tilde\cG_{lr},\tA)}|_{pt}\neq\{0\}$ iff $T_{(\tilde\cG_{lr}\tA,\tA)}|_{pt}\ssetneq T_{(\tilde\Si,\tA)}|_{pt}$. Now, by the definition of the tangent space:
\beq
T_{(\tilde\cG_{lr}\tA,\tA)}|_{pt}=\Span(U\tA+\tA V)|_{pt}+Der_\k(R_{pt})(\tA)|_{pt}\stackrel{\tA|_{pt}=\zero}{=}Der_\k(R_{pt})(\tA)|_{pt}\sseteq Mat_{(m-j)\times(n-j)}(\C).
\eeq
Thus the condition is that $Der_\k(R_{pt})(\tA)|_{pt}\ssetneq Mat_{(m-j)\times(n-j)}(\C)$. Note that the matrix structure is not used here, thus one can write the entries of $\tA$ as a column vector and form the matrix $Jac(\tA))=\{\cD_\al (\ta_{ij})\}_{(ij),\al}$, where $\{\cD_\al\}$ are some generators of $Der_\k(R_{pt})$. This matrix has $mn$ rows, while the number of columns depends on $Der_\k(R_{pt})$. Then $Der_\k(R_{pt})(\tA)|_{pt}\ssetneq Mat_{(m-j)\times(n-j)}(\C)$ iff
$I_{(m-j)(n-j)}(Jac(\tA))\ssetneq R_{pt}$. Thus
\beq
pt\in \Supp(T^1_{(\tilde\Si,\tilde\cG_{lr},\tA)})\quad \iff\quad I_1(\tA)+I_{(m-j)(n-j)}(Jac(\tA))\ssetneq R_{pt}.
\eeq
(Here we are allowed to add $I_1(\tA)$ as $\tA|_{pt}=\zero$.) Now we rewrite this condition globally. The entries of $\tA$ are precisely the local generators of $I_{j+1}(A)$ near $pt$. Thus we get:
\beq
Let\quad pt\in V\Big(I_{j+1}(A)\Big)\smin V\Big(I_j(A)\Big).\quad Then\quad
pt\in \Supp(T^1_{(\Si,\cG_{lr},A)})\quad\iff\quad
pt\in V\Big( Sing(I_{j+1}A)\Big)\smin V\Big(I_j(A)\Big).
\eeq
Together with equation \eqref{Eq.inside.proof} this finishes the proof.
\epr

\subsubsection{}{\em Proof of theorem \ref{Thm.Results.annT1.Grl.coord.changes}}

{\bf 1.} Fix some $A\in Mat_{1\times n}(R)$. The tangent space $T_{(\cG_{lr}A,A)}$ is written down in \S\ref{Sec.Background.Tangent.Spaces}. We write down the generating matrix of the submodule $T_{(\cG_{lr}A,A)}\sseteq T_{(\Si,A)}\approx Mat_{1\times n}(R)$:
\beq
\bpm
A&\zero&\dots&\dots&\zero&\{\cD_\al a_{11}\}\\
\zero&A&\zero&\dots&\zero&\{\cD_\al a_{12}\}\\
\dots&\dots\\
\zero&\dots&\dots&\zero&A&\{\cD_\al a_{1n}\}
\epm
\eeq
(The last column here denotes the block, as $\cD_\al$ run over the generators of $Der_\k(R,\cm)$.)
Then $T^1_{(\Si,\cG_{lr},A)}$ is the cokernel of this matrix, while the annihilator of $T^1_{(\Si,\cG_{lr},A)}$, i.e. the $ann.coker$ of this matrix, is precisely $Sing^{(\cm)}_n\big(I_1(A)\big)$. This proves the equality, the further embedding follows now by lemma \ref{Thm.Background.Sing(J).Properties}.

{\bf 2.} Let $\cp\ssetneq \cm$ be an ideal as in the statement. We claim that in the $\cp$-localization, $R_\cp$, the rank of $A_\cp$ is at least $(m-1)$. In fact, as $I_{m-1}(A)\not\sseteq\cp$ we get $I_{m-1}(A)_\cp\not\sseteq(\cp)_\cp$. And as $(\cp)_\cp$ is the maximal ideal of $R_\cp$ we get: $I_{m-1}(A)_\cp= R_\cp$. But then at least one of the $(m-1)\times(m-1)$ minors of $A$ becomes invertible in $R_\cp$. Therefore the matrix is equivalent to a block-diagonal,  $A_\cp\stackrel{G_{lr}}{\sim}\one_{(m-1)\times(m-1)}\oplus \tA$, with $\tA\in Mat_{1\times (n-m+1)}(R_\cp)$. Moreover, as $\cp\supseteq I_m(A)$, all the entries of $\tA$ belong to $(\cp)_\cp$.

We assume $A_\cp$ in this form, by \S\ref{Sec.Background.Invariance.Annihilator} such a transition preserves $ann(T^1_{(\Si,\cG_{lr},A)})$. Then the tangent space to the orbit (written down in \S\ref{Sec.Background.Tangent.Spaces}) decomposes into the direct sum:
\begin{multline}
(T_{(\cG_{lr}A,A)})_\cp=\Span_R\Big(UA,AV\Big)_\cp+ Der_\k(R,\cm)(A)_\cp=\Span_{R_\cp}\Big(U_\cp A_\cp,A_\cp V_\cp\Big)+Der_\k(R_\cp)(A_\cp)=\\
=Mat_{(m-1)\times n}(R_\cp)\oplus Mat_{1\times(m-1)}(R_\cp)\oplus \Span_{R_\cp}\Big(\tU\tA,\tA\tV,Der_\k(R_\cp)(\tA)\Big)_{\substack{\tU\in Mat_{1\times1}(R_\cp)\\\tV\in Mat_{(n-m+1)\times(n-m+1)}(R_\cp)}}
\end{multline}
(Here we use: $Der_\k(R_\cp)=Der_\k(R)_\cp\supseteq Der_\k(R,\cm)_\cp\supseteq(\cm\cdot Der_\k(R))_\cp=Der_\k(R_\cp)$. The last equality holds because $\cp\neq\cm$.)

We use this direct sum decomposition, together with the corresponding direct sum decomposition of $T_{(\Si,A)}$, to get:
\begin{multline}
\Big(ann(T^1_{(\Si,\cG_{lr},A)}\Big)_\cp=ann\Big((T^1_{(\Si,\cG_{lr},A)})_\cp\Big)=
ann\quotient{(T_{(\Si,A)})_\cp}{(T_{(\cG_{lr}A,A)})_\cp}\approx\\
\approx ann\quotient{Mat_{1\times (n-m+1)}(R_\cp)}{T_{(\cG_{lr}\tA,\tA)}}\stackrel{Part\ 1}{=}
Sing_{n-m+1}(I_1(\tA))
\end{multline}
Note that $I_1(\tA)=I_m(A_\cp)=I_m(A)_\cp$, and for the later ideal the expected height is $(n-m+1)$. Therefore
\beq
Sing_{n-m+1}(I_1(\tA))=Sing_{n-m+1}\big(I_m(A)_\cp\big)=Sing_{n-m+1}\big(I_m(A)\big)_\cp.
\eeq
 In the last transition we emphasize that $height(I_m(A))=height(I_m(A)_\cp)$, as $\cp\supseteq I_m(A)$.

\

{\bf 3.}  The embedding $ann(T^1_{(\Si,\cG_{r},A)})\sseteq ann(T^1_{(\Si,\cG_{lr},A)})$ is immediate as $\cG_r\sset \cG_{lr}$. The embedding
\beq
ann(T^1_{(\Si,\cG_{lr},A)})\supseteq ann.coker(A)+ann\quotient{\Mat}{Der_\k(R,\cm)(A)}
\eeq
 follows because
$\cG_{r}= G_{r}\rtimes Aut_\k(R)$ gives $T_{(\cG_{r}A,A)}\supseteq T_{(G_{r}A,A)}+T_{(Aut_\k(R)A,A)}$.

For any $0\le j<m$ we have to prove: $ann(T^1_{(\Si,\cG_{lr},A)})\sseteq Sat_{I_{j}(A)}\Big(Sing_{(m-j)(n-j)}\big(I_{j+1}(A)\big)\Big)$.

Note that $ann(T^1_{(\Si,\cG_{lr},A)})\sseteq Sat_{I_j(A)}\Big(ann(T^1_{(\Si,\cG_{lr},A)})\Big)$, thus, using part 2 of lemma \ref{Thm.Background.Localization.1}, it is enough to prove: for any prime ideal $\cp\supseteq I_{j+1}(A)$ such that $\cp\not\supseteq I_{j}(A)$ holds the inclusion $ann(T^1_{(\Si,\cG_{lr},A)})_\cp\sseteq
Sat_{I_{j}(A)}\Big(Sing_{(m-j)(n-j)}\big(I_{j+1}(A)\big)\Big)_\cp$.
 The proof below is similar to that of part two, just for $j<m-1$ we obtain weaker statements.

Take such a prime $\cp$, then the image of $I_{j}(A)$ under the $\cp$-localization is $R_\cp$, by part 0 of lemma \ref{Thm.Background.Localization.1}. Thus at least one $k\times k$ minor of $A$ is invertible in $R_\cp$. Therefore the localization is block-diagonalizable, $A_\cp\stackrel{(G_{lr})_\cp}{\sim} \one_{j\times j}\oplus\tA$, where $\tA\in Mat_{(m-j)\times(n-j)}(R_\cp)$.
By the invariance of annihilator, \S\ref{Sec.Background.Invariance.Annihilator}, we assume $A_\cp$ in this form.

Note that $I_1(\tA)=I_{j+1}(A_\cp)=I_{j+1}(A)_\cp$ and, as $I_{j+1}(A)\sseteq\cp$, $I_{j+1}(A)_\cp\sseteq(\cp)_\cp$, i.e. none of the entries of $\tA$ is invertible in $R_\cp$.
As in part two we decompose the tangent space to the orbit into the direct sum.
\begin{multline}
(T_{(\cG_{lr}A,A)})_\cp=\Span_R\Big(UA,AV\Big)_\cp+Der_\k(R,\cm)(A)_\cp=\Span_{R_\cp}\Big(U_\cp A_\cp,A_\cp V_\cp\Big)+Der_\k(R_\cp)(A_\cp)=\\
=Mat_{j\times n}(R_\cp)\oplus Mat_{(m-j)\times j}(R_\cp)\oplus \Span_{R_\cp}\Big(\tU\tA,\tA\tV,Der_\k(R_\cp)(\tA)\Big)_{\substack{\tU\in Mat_{(m-j)\times(m-j)}(R_\cp)\\\tV\in Mat_{(n-j)\times(n-j)}(R_\cp)}}.
\end{multline}
 We simplify the annihilator according to this decomposition:
\begin{multline}
\Big(ann(T^1_{(\Si,\cG_{lr},A)}\Big)_\cp=ann\Big((T^1_{(\Si,\cG_{lr},A)})_\cp\Big)=
ann\quotient{(T_{(\Si,A)})_\cp}{(T_{(\cG_{lr}A,A)})_\cp}
\approx ann\quotient{Mat_{(m-j)\times(n-j)}(R_\cp)}{T_{(\cG_{lr}\tA,\tA)}}.
\end{multline}
Unlike part two, for $j>1$ we cannot ``pack" the last term in a nice form. Instead we enlarge the annihilator ideal  by observing
that $T_{(G_{lr}\tA,\tA)}\sseteq Mat_{(m-j)\times(n-j)}(I_1(\tA))$. (This is equality for $j=1$, but can be a proper embedding for $k>1$.) Now, as in part two, we observe: $I_1(\tA)=I_{j+1}(A_\cp)=I_{j+1}(A)_\cp$. Moreover,
$height(I_{j+1}(A))=height(I_{j+1}(A)_\cp)$ and the expected height of
$I_{j+1}(A)_\cp$ is $(m-j)(n-j)$.  Therefore:
\begin{multline}
\Big(ann(T^1_{(\Si,\cG_{lr},A)}\Big)_\cp\sseteq Sing_{(m-j)(n-j)}\big(I_1(\tA)\big)=Sing_{(m-j)(n-j)}\big(I_{j+1}(A)_\cp\big)=
\Big(Sing_{(m-j)(n-j)}\big(I_{j+1}(A)\big)\Big)_\cp.
\end{multline}
As this embedding holds for any localization at $\cp\not\supseteq I_j(A)$, we get:
 $ann(T^1_{(\Si,\cG_{lr},A)})\sseteq Sat_{I_j(A)}\Big(Sing_{(m-j)(n-j)}\big(I_{j+1}(A)\big)\Big)$.

\

{\bf 3'.} Follows right from part 4 of lemma \ref{Thm.Background.Sing(J).Properties} and part 3.

\

{\bf 4.}
 The embedding $\sseteq$ follows from part three. For the embedding $\supseteq$ we use
 lemma \ref{Thm.Background.Localization.2}. Thus it is enough to verify that
 for any non-maximal prime ideal, $\cp\ssetneq\cm$, the localizations of the ideals  satisfy:
\beq
\text{If}\quad
ann(T^1_{(\Si,\cG_{lr},A)})_\cp\neq R_\cp\quad \text{then} \quad
\capl^{m-1}_{r=0}Sat_{I_r(A)}\Big( Sing_{(m-r)(n-r)}(I_{r+1}(A))\Big)_\cp\neq R_\cp.
\eeq

 Fix a prime ideal $\cp\ssetneq\cm$ and fix the number $j$ satisfying $\cp\supseteq I_{j+1}(A)$, $\cp\not\supseteq I_j(A)$. It exists, because $I_0(A)=R$, $I_{m+1}(A)=\{0\}$, and moreover $0\leq j\le m$. Furthermore, as the chain of ideals $I_r(A)$ is monotonic, $j$ is uniquely defined by $\cp$.

Note that $I_{r+1}(A)_\cp=R_\cp$ for $r+1\le j$, because $I_{r+1}(A)\not\sseteq\cp$, lemma \ref{Thm.Background.Localization.1}. Therefore
$Sat_{I_r(A)}\Big( Sing_{(m-r)(n-r)}(I_{r+1}(A))\Big)_\cp=R_\cp$ for $r+1\le j$.
Thus for the localization at $\cp$ it is enough to prove:
\beq
\text{If}\quad
ann(T^1_{(\Si,\cG_{lr},A)})_\cp\neq R_\cp\quad \text{then} \quad
Sat_{I_j(A)}\Big( Sing_{(m-j)(n-j)}(I_{j+1}(A))\Big)_\cp\neq R_\cp.
\eeq

As  $\cp\not\supseteq I_j(A)$ we have $I_j(A)_\cp=R_\cp$, thus the localization $A_\cp$ of $A$ has at least one invertible minor of size $j\times j$. Thus  $A_\cp$ is  $(G_{lr})_\cp$-equivalent to $\one_{j\times j}\oplus\tA$, where $\tA\in Mat_{(m-j)\times(n-j)}(\cp_\cp)$. Recall that $ann(T^1_{(\Si,\cG_{lr},A)})$ is invariant under the $G_{lr}$-equivalence, \S\ref{Sec.Background.Invariance.Annihilator}, therefore from now on we assume $A_\cp$ in this form.

For this form of $A_\cp$ we have the direct sum decomposition
\beq
T_{((\cG_{lr})_\cp,A_\cp)}\approx Mat_{j\times n}(R_\cp)\oplus Mat_{(m-j)\times j}(R_\cp)\oplus T_{(\widetilde{\cG_{lr}},\tA)},
\eeq
 as in part 2. Thus
 $ann(T^1_{(\Si,\cG_{lr},A)})_\cp\approx ann(T^1_{(\tilde\Si,\tilde\cG_{lr},\tA)})$, where $\tilde\Si=Mat_{(m-j)\times(n-j)}(R_\cp)$ and $\tilde\cG_{lr}$ is the corresponding group. Therefore we must prove:
\beq
\text{If}\quad
ann(T^1_{(\tilde\Si,\tilde\cG_{lr},\tA)})\neq R_\cp\quad \text{then}\quad
 Sat_{I_j(A)}\Big( Sing_{(m-j)(n-j)}(I_{j+1}(A))\Big)_\cp\neq R_\cp.
\eeq
Recall that $T_{(\widetilde{\cG_{lr}},\tA)}=\Span_{R_\cp}(U\tA+\tA V)+Der_\k(R_\cp)(\tA)$ and all the entries of $\tA$ are in $\cp$. Therefore $T_{(\widetilde{\cG_{lr}},\tA)}\ssetneq T_{(\tilde\Si,\tA)}$ iff
$I_1(\tA)\cdot T_{(\tilde\Si,\tA)}+Der_\k(R_\cp)(\tA)\ssetneq T_{(\tilde\Si,\tA)}$. Thus $ann(T^1_{(\tilde\Si,\tilde\cG_{lr},\tA)})\neq R_\cp$ iff $Sing(I_{1}(\tA))\neq R_\cp$. Finally note that the expected height of $I_1(\tA)$ is precisely ${(m-j)(n-j)}$ and $I_1(\tA)=I_{j+1}(A)_\cp$, thus
$Sing(I_1(\tA))=Sing(I_{j+1}(A))_\cp$.
\epr

\bex\label{Ex.Results.Computation.sqrt(ann)}
Let $R=\k[\![x_1,\dots,x_p]\!]$, let $A=\bpm x_2&x^k_1\\x^l_1&x_2\epm$. By part 4 of theorem \ref{Thm.Results.annT1.Grl.coord.changes} we have:
\beq
ann(T^1_{(\Si,\cG_{lr},A)})\sseteq
\Big(Sat_{I_1(A)}\big(Sing_1^{(\cm)}(I_2(A))\big)\cap Sat_{I_0(A)}\big(Sing_4^{(\cm)}(I_1(A))\big)\Big)
 \sseteq\sqrt{ann(T^1_{(\Si,\cG_{lr},A)})}
\eeq
 For $Sat_{I_1(A)}\big(Sing_1^{(\cm)}(I_2(A))\big)$ we have:
\begin{multline}
I_2(A)=(det(A))=(x^2_2-x^{k+l}_1),\quad Jac^{(\cm)}(I_2(A))=\underbrace{\big((k+l)x^{k+l}_1,(k+l)x^{k+l-1}_1x_2,x_1x_2,x^2_2,0,\dots,0\big)}_{p},\quad
\\Der_\k(R,\cm)(I_2(A))=\cm\cdot\big((x^{k+l-1}_1)+(x_2)\big)\sset R,\quad
Sing_1(I_2(A))=\cm\cdot\big((x^{k+l-1}_1)+(x_2)\big),
\end{multline}
and therefore $Sat_{I_1(A)}\Big(Sing_1(I_2(A))\Big)=R$.

To compute $Sat_{I_0(A)}\big(Sing_4^{(\cm)}(I_1(A))\big)$ we note that $rank(Der_\k(R))=dim(R)<4$, therefore
$Sing_4^{(\cm)}(I_1(A))=I_1(A)=(x_2)+(x^{k}_1)+(x^{l}_1)$. Finally $I_0(A)=R$, thus
$Sat_{I_0(A)}\big(Sing_4^{(\cm)}(I_1(A))\big)=(x_2)+(x^{k}_1)+(x^{l}_1)$. Altogether we get:
\beq
ann(T^1_{(\Si,\cG_{lr},A)})\sseteq
(x_2)+(x^{k}_1)+(x^{l}_1)
 \sseteq\sqrt{ann(T^1_{(\Si,\cG_{lr},A)})}=(x_2)+(x_1)
\eeq
 In particular,  $A$ is finitely-$\cG_{lr}$-determined iff $\sqrt{ann\big(T^1_{(\Si,\cG_{lr}A,A)}\big)}=\cm$ iff $p=2$.
\eex
\bex
More generally, let $R=\k[\![x_1,\dots,x_p]\!]$ and suppose $I_1(A)\sseteq(x_1,\dots,x_{p-1})^2$.
 Then $Der_\k(R)(I_1(A))\sseteq (x_1,\dots,x_{p-1})R^{\oplus k}$, giving $Sing(I_1(A))\sseteq (x_1,\dots,x_{p-1})$.
 Then theorem \ref{Thm.Results.annT1.Grl.coord.changes} gives (for $j=1$): $ann\big(T^1_{(\Si,\cG_{lr},A)}\big)\sseteq(x_1,\dots,x_{p-1})\sset R$. In particular this ideal  cannot contain any power of $\cm$. Therefore $A$ is not finitely-$\cG_{lr}$-determined.

Even more generally, for any Noetherian ring $R$, suppose $I_1(A)\sseteq J^2$, where $\sqrt{J}\ssetneq\cm$. Then $Sing(I_1(A))\sseteq J$, thus $ann\big(T^1_{(\Si,\cG_{lr},A)}\big)\sseteq J$, thus $ann\big(T^1_{(\Si,\cG_{lr},A)}\big)$ cannot contain any power of maximal ideal.

Yet more generally, suppose for some $0\le j<m$ holds: $I_{j+1}(A)\sseteq J^2$, with $\sqrt{J}\ssetneq\cm$. Suppose in addition that $height(I_{j+1}(A))<height(I_j(A))$. Then $A$ is not finitely determined.
\eex

\subsection{Finite determinacy for $\cG_{lr}$}

\

{\em \underline{Proof of proposition \ref{Thm.Results.Finite.Determinacy.Glr.with.coord.change.Corol}}}

{\bf 1.}
If $rank(Der_\k(R))<(n-j)(m-j)$ then, by part 3'  of theorem \ref{Thm.Results.annT1.Grl.coord.changes},
 $ann T^1_{(\Si,\cG_{lr},A)}\sseteq Sat_{I_j(A)}(I_{j+1}(A))$.
 But $height(I_{j+1}(A))\le (n-j)(m-j)<dim(R)$, as $A\in\Matm$, thus $ann( T^1_{(\Si,\cG_{lr},A)})$ cannot contain any power of $\cm$. Thus $A$ is not finitely-$\cG_{lr}$-determined.

 \

{\bf 2.}
Recall that the determinantal ideals are preserved under the action of $G_{lr}$. The group $Aut_\k(R)$ does not preserve them but preserves their heights. Therefore, the statement follows straight from proposition \ref{Thm.background.Fitt.Ideals.generic.height}.


\

{\bf 3.}
As $G_r\sseteq\cG_{lr}$, the only non-trivial direction is that $\cG_{lr}$-determinacy implies the $G_{r}$-determinacy. If $A$ is $\cG_{lr}$-finitely determined then the ideal of maximal minors is of the expected height, $height(I_m(A))=min\Big(n-m+1,dim(R)\Big)$. As $dim(R)\le n-m+1$ we get $height(I_m(A))=dim(R)$,
 thus $I_m(A)$ contains a power of the maximal ideal. But then, by proposition 2.3 of \cite{Belitski-Kerner1}, $A$ is finitely $G_r$-determined.

\

{\bf 4.} By the assumption: $B\in \sqrt{J}\cdot T_{(\cG_{lr}A,A)}\sseteq T_{(\cG^{(\sqrt{J})}_{lr}A,A)}$.
 Therefore the statement follows from proposision \ref{Thm.Background.Reduction.to.annihilator}. (We just note that the filtration $Mat_{m\times n}((\sqrt{J})^q)$ is $\cG^{(\sqrt{J})}_{lr}$-invariant and moreover, the group $\cG^{(\sqrt{J})}_{lr}$ is unipotent for this filtration.)
\epr

{\em\underline{Proof of proposition \ref{Thm.Intro.Glr.coord.changes.Additional.Corollary}}}

1. We should check whether the ideal $ann(T^1_{(\Si,\cG_{lr},A)})$ contains some finite power of maximal ideal or geometrically whether $Supp(T^1_{(\Si,\cG_{lr},A)})$ is the origin. For this we use proposition \ref{Thm.Results.Support.Set.Theoretic.Glr.coord.change}.

If $A$ is finitely determined and $\k=\bar{\k}$ the dimension of $V(I_{r}(A))$ is the expected one.
As $dim(R)\le 2(n-m+2)$ this means: $V(I_r(A))$ is one point for $r<m$. Therefore of all the terms in the right hand side of the equality in proposition \ref{Thm.Results.Support.Set.Theoretic.Glr.coord.change} we should check only the term $\overline{Sing(V_m(A))\smin V(I_{m-1}(A))}$. Thus the support of $ann(T^1_{(\Si,\cG_{lr},A)})$ is a point iff
$Sing(V_m(A))\smin V(I_{m-1}(A))$. Hence the statement.

2. This follows straight from part one, note that $dim(V(I_m(A)))=1$. We should only add: as $Spec(R)$ is smooth, and the determinantal ideals $I_j(A)$ are of expected
height, their zero loci $V(I_j(A))$ are Cohen-Macaulay. In particular they have no embedded components.

3. Apply proposition \ref{Thm.Results.Support.Set.Theoretic.Glr.coord.change} for the $j=1$ case.
\epr

\bex\label{Example.Non.CM.ring}
If $R$ is a regular ring and $A\in \Mat$ is finitely determined, then, in particular,  the ideals $\{I_{j+1}(A)\}$ for which  $height(I_j(A))> height(I_{j+1}(A))$ are radical.
\li If $R$ is non-regular then being finitely determined does not imply that the ideal $I_1(A)$ is radical
 (i.e. the zero locus $A^{-1}(0)$ is reduced). For example, let $R=\quotients{\k[\![x,y,z]\!]}{(xz,yz)}$ and $A=x+y+z\in Mat(1,1,R)$.
 Then $A$ is obviously $Aut_\k(R)$-finitely-determined (even stable), but $I_1(A)$ defines a non-reduced scheme, whose local ring
is $\quotients{\k[\![x,y,z]\!]}{(xz,yz,x+y+z)}\approx\quotients{\k[\![x,y]\!]}{(x^2,xy)}$.
\li Further, see remark \ref{Ex.Fin.Det.Nonisolated.Sing.for.non-regular.ring}, even if $R$ is a complete
intersection and $A$ is finitely-$\cG_{lr}$-determined, the ideal $I_1(A)$ can define a scheme
with multiple components.
\eex

\bprop\label{Thm.Intro.Fin.Det.vs.Stability.Nearby}
Let $R=\quotients{\C\{\ux\}}{I}$, then
$A\in\Mat$ is finitely-$\cG_{lr}$-determined iff for each point $0\neq pt\in Spec(R)$
 the matrix $A$ is ${\cG_{lr}}$-stable at $pt$.
\eprop
(In this statement by $Spec(R)$ we mean a small enough neighborhood of the origin, so that $A$ is defined at each point of $Spec(R)$. By the group action of ${\cG_{lr}}$ at $pt$ we mean $GL(m,R_{pt})\times GL(n,R_{pt})\rtimes \cR_{pt}$, where $R_{pt}$ is the local ring, while $\cR_{pt}$ is the group of local coordinate changes that do not necessarily preserve $pt$, i.e. $\cR$ includes the translations.)

{\em\underline{Proof of corollary \ref{Thm.Intro.Fin.Det.vs.Stability.Nearby}}.}
As is mentioned in the introduction, the finite determinacy is equivalent to $\cm^N\sseteq ann(T^1_{(\Si,G,A)})$.
 Thus $A$ is finitely-$\cG_{lr}$-determined iff the module $T^1_{(\Si,G,A)}$ is supported at the origin only, i.e. for any $pt\in Spec(R)\smin\{\o\}$ holds: $T_{(\Si,A)}|_{pt}=T_{(\cG_{lr},A)}|_{pt}$, which is precisely the local stability.
\epr

\subsection{The congruence action, $\cG_{congr}$}

\subsubsection{The support of $T^1_{(\Si,\cG_{congr},A)}$}  (The proof of theorem \ref{Thm.Results.T1.for.Congruence}.)

\bpr
{\bf 1.} It is enough to prove: if $rank(Der_\k(R))<\lfloor\frac{m}{2}\rfloor$ then the module $T_{(\cG_{congr}A,A)}$ is of rank smaller than $rank(T_{(\Si,A)})=m^2$. And thus the module $T^1_{(\Si,\cG_{congr},A)}$ is of positive rank, in particular is generically supported on $Spec(R)$, i.e. $ann(T^1_{(\Si,\cG_{congr},A)})=\{0\}$.

To compute the rank we localize at the generic point of $Spec(R)$. More precisely, choose an irreducible component and localize at the ideal of all the nilpotent elements. Then $R_{(0)}$ is a field, the modules $T_{(\cG_{congr}A,A)}$, $T_{(\Si,A)}$ become  vector spaces and the relevant ranks are their dimensions.

To give an upper bound on the dimension of $(T_{(\cG_{congr}A,A)})_{(0)}=\Span_{R_{(0)}}(UA+AU^T)+Der_\k(R)(A)_{(0)}$ we need the lower bound on the dimension of the vector space
\beq
\{U\in Mat_{m\times m}(R_{(0)})|\ UA+AU^T=\zero\}.
\eeq
The later equation is well studied, the dimension of the space of solutions is precisely the codimension of the orbit of $A$ under the congruence. The minimal codimension equals
$\lfloor\frac{m}{2}\rfloor$, see e.g. Theorem 3 in \cite{De Teran-Dopico}, and it is achieved for $A\in Mat_{m\times m}(R_{(0)})$ generic. Therefore we get: $rank(T_{(G_{congr}A,A)})\le m- \lfloor\frac{m}{2}\rfloor$. Hence the stated bound.

\

The proofs of the remaining parts are essentially the same as in theorem \ref{Thm.Results.annT1.Grl.coord.changes}, thus we give only the sketches.

{\bf 2.i}
Localize at $\cp$, then bring $A_\cp$ to the canonical form $A\stackrel{G_{congr}}{\sim}\one_{(m-1)\times(m-1)}\oplus \tA_{1\times 1}$, see  \cite[\S3.1]{Belitski-Kerner1}. Now, as in the proof of theorem \ref{Thm.Results.annT1.Grl.coord.changes}, one has $(T^1_{(\Si,\cG_{congr},A)})_\cp\approx\quotient{R_\cp}{(\tA)+(Jac(\tA))}=\quotient{R_\cp}{(det(A))_\cp+(Jac(det(A))_\cp)}$.
Thus $ann(T^1_{(\Si,\cG_{congr},A)})_\cp=Sing_1(det(A))_\cp$.

For parts 2.ii, 2.iii we note that the expected height of $I_{j+1}(A)$ is $min\Big(\bin{m-j+1}{2},dim(R)\Big)$.

{\bf 2.ii} The same as in theorem \ref{Thm.Results.annT1.Grl.coord.changes}, just replace $G_{lr}$ by $G_{congr}$.

{\bf 2.iii} Follows straight from part 4 of lemma \ref{Thm.Background.Sing(J).Properties}.

{\bf 3.i} In this case $A_\cp\stackrel{G_{congr}}{\sim}E_{(m-2)\times(m-2)}\oplus\tA_{2\times 2}$, where $E_{(m-2)\times(m-2)}$ is skew-symmetric and invertible, while $\tA\in Mat^{skew-sym}_{2\times 2}(\cp_\cp)$,
see \cite[\S3.1]{Belitski-Kerner1}. Then as above we get: $ann(T^1_{(\Si,\cG_{congr},A)})_\cp=Sing_1(I_1(\tA))=Sing_1(I_{m-1}(A))$.

{\bf 3.ii} As in the proof of theorem \ref{Thm.Results.annT1.Grl.coord.changes} we fix a prime ideal
$I_j(A)\not\sseteq\cp\supseteq I_{j+1}(A)$. Note that for $j$ odd one has $\sqrt{I_j(A)}=\sqrt{I_{j+1}(A)}$, thus $\cp$ as above exists only for $j$-even. Otherwise the proof is the same.
\epr

\bcor
Suppose  $\k=\bk$ and $R$ is a regular ring, Noetherian and Henselian. Then $A\in Mat^{sym}_{m\times m}(R)$ is finitely-$(\Si^{sym},\cG_{congr})$-determined  iff all the degeneracy loci $V(I_j(A))\sset Spec(R)$ are of expected dimensions (or empty) and all the complements $V(I_j(A))\smin V(I_{j-1}(A))$ are smooth.
\ecor
(The criterion for skew-symmetric matrices is formulated similarly.)

\end{document}